\documentclass[12pt,thmsa]{article}
\usepackage{amsfonts}
\usepackage[dvips]{graphics}

\newtheorem{theorem}{Theorem}[section]
\newtheorem{lemma}[theorem]{Lemma}
\newtheorem{corollary}[theorem]{Corollary}
\newtheorem{proposition}[theorem]{Proposition}
\newtheorem{example}[theorem]{Example}
\newtheorem{remark}[theorem]{Remark}

\newtheorem{hypothesis}[theorem]{Hypothesis}

\def\bit{\begin{itemize}}
\def\eit{\end{itemize}}
\reversemarginpar   
\def\bc{\begin{center}}
\def\ec{\end{center}}
\def\bthm{\begin{theorem}}
\def\ethm{\end{theorem}}
\def\bcor{\begin{corollary}}
\def\ecor{\end{corollary}}
\def\bprop{\begin{proposition}}
\def\eprop{\end{proposition}}
\def\blem{\begin{lemma}}
\def\elem{\end{lemma}}
\def\bex{\begin{example} {\rm }
\def\eex{\end{example} }}
\def\brem{\begin{remark}}
\def\erem{\end{remark}}
\def\prf{\noindent{\bf Proof~: }}
\def\bdes{\begin{description}}
\def\edes{\end{description}}
\def\ita{\item[(a)]}
\def\itb{\item[(b)]}
\def\itc{\item[(c)]}
\def\itd{\item[(d)]}

\def\iti{\item[(i)]}
\def\itii{\item[(ii)]}
\def\itiii{\item[(iii)]}
\def\itiv{\item[(iv)]}
\def\itv{\item[(v)]}
\def\itvi{\item[(vi)]}

\def\beq{\begin{equation}}
\def\eeq{\end{equation}}

\def\ben{\begin{enumerate}}
\def\een{\end{enumerate}}

\def\beqar{\begin{eqnarray}}
\def\eeqar{\end{eqnarray}}
\def\beqarr{\begin{eqnarray*}}
\def\eeqarr{\end{eqnarray*}}
\def \non{{\nonumber}}

\def\RR{{\mathbb R}}  

 \def\cB{\mathcal{B}}    
\def\cD{\mathcal{D}}    
 \def\cH{\mathcal{H}}   
    
\def\cM{\mathcal{M}} \def\cN{\mathcal{N}}   
\def\cP{\mathcal{P}}

\def\qed{\hspace{.1in}{\bf QED}}

\def\P{{\mathsf P}} 

\def\E{{\mathsf E}} 

\def\ZZ{{\mathbb Z}}       

\def\one{{\bf 1}}

\def\rar{\rightarrow}

\def\vareps{\varepsilon}
\def\eps{\epsilon}

\def\la{\langle}

\def\ra{\rangle}

\def\dist{\hbox{dist}}

\def\part{\partial}

\def\d#1dt{\frac{d#1}{dt}}    

\begin{document}
\title{Self Interacting Diffusions \\ III: Symmetric Interactions
\author{{\bf Michel Bena\"{\i}m}\\Institut de Math\'ematiques\\
Universit\'e de Neuch\^atel, Suisse \and
{\bf Olivier Raimond}\\Laboratoire de  Mod\'elisation Stochastique
et Statistique\\Universit\'e Paris Sud, France}}
\date{}
\maketitle
\bibliographystyle{apalike}
\begin{center} Dedicated to Morris W Hirsch 70's birthday \end{center}
\begin{abstract}
Let $M$ be a compact Riemannian manifold. A {\em self-interacting
  diffusion} on $M$ is a stochastic process solution to 
$$dX_t = dW_t(X_t) - \frac{1}{t}\left(\int_0^t \nabla V_{X_s}(X_t)ds\right)dt$$
where $\{W_t\}$ is a Brownian vector field on $M$ and $V_x(y) =
V(x,y)$ a smooth function. Let $\mu_t = \frac{1}{t} \int_0^t
\delta_{X_s} ds$ denote  the normalized occupation measure of
$X_t$. We prove that, when $V$ is symmetric, $\mu_t$ converges almost
surely to the critical set of a certain nonlinear free energy
functional $J$. Furthermore,  $J$ has generically finitely many
critical points and $\mu_t$ converges almost surely toward a local
minimum of $J.$ Each local minimum having a positive probability to be
selected.
\end{abstract}
{\bf Acknowledgment:}
We are very grateful to Gerard Ben Arous, Thierry Coulhon, Morris W
Hirsch, Josef Hofbauer, Florent Malrieu and  Hans Henrik Rugh for
their suggestions and comments.
\tableofcontents

\section{Introduction}
\label{intro}
Let $M$ be a $C^{\infty}$ $d$-dimensional, compact connected 
 Riemannian manifold without boundary and $V : M \times M \to \RR$
be a  smooth function called a  {\em potential.}
For every Borel probability measure $\mu$ on $M$ let 
$V{\mu} : M \to \RR$ denote  the smooth function defined by
\beq V{\mu}(x) = \int_M V(x,u)\mu(du),\eeq
and let $\nabla (V{\mu})$ denote its gradient (computed with respect to
the Riemannian metric on $M$).

A {\em Self-interacting diffusion process} associated to $V$  is a
continuous time stochastic process  living on $M$  solution  to the
stochastic differential equation (SDE)
\begin{equation}
\label{eq:sde}
dX_t = \sum_{i = 1}^N F_i(X_t) \circ dB_t^i - \frac{1}{2} \nabla
(V{\mu_t})(X_t) dt,\quad X_0 = x \in M
\end{equation}
where $(B^1, \ldots, B^N)$ is a standard Brownian motion on $\RR^N$,
$\{F_i\}$ is a family of smooth vector fields on $M$ such that
\beq \sum_{i=1}^N F_{i} (F_i f) = \Delta f \eeq
(for $f \in C^{\infty}(M)$), where $\Delta$ denotes the Laplacian on
$M$; and
\beq \mu_t = \frac{1}{t}\int_0^t \delta_{X_s} ds \eeq
is  the {\em empirical occupation measure} of $\{X_t\}$.

In absence of drift (i.e~ $V(x,y) = 0$) $\{X_t\}$  is just a Brownian motion on $M.$ If $V(x,y) = V(x)$ then it is a diffusion process on $M.$ However, for a general function $V,$ such a process 
is  characterized by the fact that the drift term in equation
(\ref{eq:sde}) depends both on the position of the process  and
its empirical occupation measure up to time $t.$ 

Self-interacting diffusions (as defined here) were introduced
in  Benaim, Ledoux and Raimond (2000), (hereafter referred  as (BLR))  and we  refer the reader to this paper for a  more detailed
definition and  basic properties. 

It is worth pointing out that equation (\ref{eq:sde}) presents some strong similarities  with the following class of SDE
\beq
\label{eq:bp}
dY_t = dB_t - \left ( \int_0^t v'(Y_s-Y_t)ds \right )dt
\eeq
 whose behavior has been the focus of much attention in the recent years  
 (see e.g~Norris, Williams and Rogers (1987),  Durret and Rogers
 (1992), Cranston and LeJan (1995), Cranston and Mountford (1996),
 Raimond (1997), Hermann and Roynette (2003) or Pemantle (2002)
 for a recent overview and further references about reinforced random processes). The main differences
 being that
\bdes
\iti The SDE (\ref{eq:sde})
lives on an arbitrary but compact manifold, while (\ref{eq:bp}) lives on $\RR$ or $\RR^d.$
\itii The  drift term in (\ref{eq:bp})  depends on the {\bf non-normalized} occupation measure
 $$t \mu_t = \int_0^t  \delta_{X_s} ds.$$
\edes
A major goal in understanding (\ref{eq:sde}) is
\bdes
 \ita to provide  tools allowing to  analyze the long term behavior of
$\{\mu_t\};$ and, using these tools,
\itb  to identify (at least partially) general classes of potential leading to certain types of behaviors.
\edes
A first step in this direction has been achieved in
(BLR), where it is shown that the asymptotic behavior of $\{\mu_t\}$
can be precisely described in terms of  a certain deterministic semi-flow
$\Psi = \{\Psi_t\}_{t \geq 0}$  defined on the space of Borel
probability measures on $M.$ 
For instance, there are situations (depending on the shape of $V$) in
which $\{\mu_t\}$ converges almost surely to an equilibrium point
$\mu_\infty$  of $\Psi$ ($\mu_\infty$ is random) and other situations where the limit set of
$\{\mu_t\}$ coincides almost surely with a periodic orbit for $\Psi$
(see the examples in section 4 of (BLR)). 

The present paper adresses the second part of this program.
The main result here is  that 
\begin{quote} {\em Symmetric interactions (i.e~ symmetric potentials)  force  $\{\mu_t\}$   to converge  almost surely toward the critical set of a certain nonlinear free-energy functional}. 
\end{quote} 
This result  encompasses most of the examples considered
in (BLR) and   enlightens the results of (BLR) and Benaim and Raimond (2002).
It also  
allows to give a sensible definition of {\em self-attracting} or {\em
repelling} diffusions.  

The organization of the paper is as follows.
Section \ref{sec:hypo} defines the class of  potentials considered here, gives some examples and states the main results.  Section \ref{sec:background} reviews some material  from (BLR) on which rely the analysis.  Sections \ref{sec:proof1}, \ref{sec:attract}, \ref{sec:nonconv} and \ref{sec:appendix} are devoted to the proofs.
\section{Hypotheses and main results}
\label{sec:hypo}
We assume throughout that  $V$ is a $C^3$ map\footnote{This regularity condition can be slightly
weakened (see Hypothesis 1.4  in (BLR)).} 
and that 
\begin{hypothesis} [Standing assumption]
\label{hyp:main}
$V$ is   symmetric :  $$V(x,y) = V(y,x).$$
\end{hypothesis}

Recall that $\lambda$ denotes the Riemannian probability on $M.$
We will sometime use the following additional hypothesis:
\begin{hypothesis} [Occasional  assumption 1]
\label{hyp:occ}
The mapping
\beq V\lambda : x \mapsto V\lambda(x) = \int_M V(x,y)\lambda(dy)\eeq
is constant.
\end{hypothesis}

This later condition has the interpretation that if the empirical
occupation measure of $X_t$ is (close to) $\lambda$ then the drift
term $\nabla V \mu_t (X_t)$ is (close to) zero. In other words, if the
process has visited $M$ ``uniformly'' between times $0$ and $t$ then
it has no preferred directions and behaves like a Brownian
motion.

\paragraph{Notation.} 
Throughout we let $C^0(M)$ denote the Banach space of real valued continuous functions $f  : M \to \RR,$ equipped with the  supremum norm
$$||f||_{\infty} = \sup_{x \in M} |f(x)|.$$
Given a positive function $g \in C^0(M)$
we let $\la \cdot  , \cdot \ra_{g}$ denote the inner product on $C^0(M)$
defined by
$$
\la u,v \ra_{g} = \int_M u(x)v(x)g(x)\lambda(dx).$$
When $g = 1$ we usually write
$\la \cdot , \cdot \ra_{\lambda}$ (instead of $\la \cdot ,\cdot  \ra_{1}$) and 
$||f||_{\lambda}$ for $\sqrt{\la f, f \ra_{\lambda}}.$ 

The completion of $C^0(M)$ for the norm $||f||_{\lambda}$ 
is the Hilbert space $L^2(\lambda).$
We sometime use the notation $\one$ to denote the function on $M$ taking value one everywhere; and
 $$L^2_0(\lambda) = \one^{\perp} = \{h \in L^2(\lambda) \: : \la h, \one \ra_{\lambda} = 0\}.$$
We let  ${\cal M}(M)$ denote the space of Borel
bounded measures on $M$ and ${\cal P}(M)$ the subset of Borel probabilities.  For $\mu \in {\cal M}(M)$ and
$ f \in C^0(M)$ we set 
\beq \mu f = \int_M f(x) \mu(dx)\eeq
and 
\beq |\mu| = \sup \{ |\mu f|: \: f \in C^0(M), ~\|f\|_{\infty} =
1\}. \eeq
We let  ${\cal M}_s(M)$ denote the Banach space $({\cal M}(M),
|\cdot|)$   (i.e., the dual of  $C^0(M)$) and
${\cal M}_w(M)$ (respectively,  ${\cal P}_w(M)$)   the metric space obtained by equipping ${\cal
M}(M)$ (respectively,  ${\cal P}(M)$) with the narrow (or weak*) topology. In particular, ${\cal P}_w(M)$ is a compact subspace of ${\cal M}_w(M).$
Recall that the narrow topology is the topology 
induced by the family of semi-norms
$\{\mu \mapsto |\mu f|: \, f \in C^0(M)\}.$
Hence $\mu_n \to \mu$ in ${\cal M}_w(M)$ if and only if $\mu_n f \to \mu f$
for all $f \in C^0(M).$

Everywhere in the paper  a subset of a topological space inherits the induced topology.

\paragraph{The operator $V.$}
The function $V$  induces an operator
$$ V : {\cal M}_s(M) \to C^0(M),$$
defined by
\beq
\label{eq:defV}
V\mu(x) = \int_M V(x,y) \mu(dy).
\eeq
If $g \in L^2(\lambda)$ we write $Vg$ for $V(g\lambda),$  where $g
\lambda$ stands for the measure whose Radon Nikodym derivative with
respect to $\lambda$ is $g$.

\medskip
The following basic lemma will be used in several places
\blem \label{th:Vcompact}
\bdes
\iti
The operator  $V : {\cal M}_s(M) \to C^0(M)$ 
and its restriction to $L^2(\lambda)$ (defined by $g \mapsto
V(g\lambda)$) are compact operators.
\itii $V$ maps continously ${\cal P}_w(M)$ into $C^0(M).$
\edes
\elem
\prf
$(i)$ Let $\mu \in {\cal M}_s(M)$. 
Then $\|V\mu\|_{\infty} \leq \|V\|_{\infty} |\mu|$ and
$|V\mu(u) - V\mu(v))| \leq (\sup_{z \in M} |V(u,z) - V(v,z)|)|\mu|$.
Therefore the set $\{ V\mu \, : |\mu| \leq 1\}$ is bounded and
equicontinuous, hence, relatively compact in $C^0(M)$ by Ascoli's
theorem. 
This proves that $V$ is compact. 

By definition $V|L^2(\lambda)$ is the composition of $V$ with the
bounded operator $g \in L^2(\lambda) \to g \lambda \in {\cal M}_s(M).$
It is then compact. 

$(ii)$ Let $\{\mu_n\}$ be a converging sequence in ${\cal P}_w(M)$ and $\mu =
\lim_{n \to \infty} \mu_n$.
Narrow convergence implies that $V\mu_n(u) \to V\mu(u)$ for all $u \in
M.$ Since, by $(i)$, $\{V\mu_n\}$ is relatively
compact in $C^0(M),$ it follows that $V\mu_n \to V\mu$ in $C^0(M)$.
\qed
\subsection{The global convergence theorem}
Let $\Pi = \Pi_V : {\cal P}_w(M) \to {\cal P}_w(M)$ be the map\footnote{
We  use the notation $\Pi_V$ for $\Pi$ when we want to emphasize the dependency on $V.$} defined by
\beq \label{eq:PI}
\Pi(\mu)(dx)  = \xi( V{\mu})(x) \lambda(dx) 
\eeq
where $\xi : C^0(M) \to C^0(M)$ is the function defined by
\beq \label{eq:defxi}
\xi(f)(x) = \frac{e^{-f(x)}}{\int_M e^{-f(y)}\lambda(dy)}.
\eeq
The {\em limit set} of $\{\mu_t\}$ denoted $L(\{\mu_t\})$ is the set
of limits (in ${\cal P}_w(M)$) of convergent sequences
$\{\mu_{t_k}\}$, $t_k \to \infty$.

\smallskip
The following theorem describes $L(\{\mu_t\})$ in terms of  $\Pi.$
 It is proved in section \ref{sec:proof1}.
\bthm \label{th:main}
With probability one $L(\{\mu_t\})$ is a compact connected subset of 
\beq {\mathsf {Fix}}(\Pi) = \{\mu \in {\cal P}_w(M): \: \mu =
\Pi(\mu)\}.\eeq
\ethm

This clearly implies 
\bcor  Assume $\Pi$ has isolated fixed points. Then  $\{\mu_t\}$ converges almost surely to a fixed point of $\Pi.$ 
\ecor 
\brem {\em  By Theorem \ref{th:generic} below, $\Pi$ has  {\em generically} isolated fixed points. Hence, the generic behavior of $\{\mu_t\}$ is convergence toward one of those fixed points.}
\erem
\subsection{Fixed points of $\Pi$}
With Theorem \ref{th:main} in hands, it is  clear  that our
description of self-interacting diffusions (satisfying hypothesis
\ref{hyp:main}) on $M$   relies on  our understanding of  the fixed 
points structure of $\Pi$. 

Let $${\cal B}_1 = \{f \in C^0(M) : \: \la f,\one \ra_{\lambda} = 1\}$$
and
$${\cal B}_0 = \{f \in C^0(M) : \: \la f,\one \ra_{\lambda} =  0\}.$$
Spaces ${\cal B}_0$ and ${\cal B}_1$ are respectively a Banach space and 
a Banach affine space parallel to  ${\cal B}_0.$

Let $$X = X_V : {\cal B}_1 \to {\cal B}_0$$
be the $C^{\infty}$ vector field defined by
\beq
\label{defX}
X(f) = - f + \xi(Vf).
\eeq
The following lemma relates fixed points of $\Pi$ to the zeroes of $X$.
\blem
\label{th:fixzero}
 Let $\mu \in {\cal P}(M).$ Then, $\mu$ is a fixed point of $\Pi$ if and only if $\mu$ is absolutely continuous with respect to $\lambda$ and  $\frac{d\mu}{d\lambda}$ is a zero of $X.$ 
 Furthermore, the map 
\beq \begin{array}{lllll}
j &:&  \mathsf{Fix}(\Pi) &\to& X^{-1}(0)\\
&& \mu  &\mapsto& \frac{d\mu}{d\lambda} \end{array}\eeq
is an homeomorphism.
In particular, $X^{-1}(0)$ is compact.
\elem
\prf
The first assertion is  immediate from the definitions. Continuity of $j$ follows from the  continuity of $\xi$ and  Lemma \ref{th:Vcompact}, $(ii).$ 
Continuity of $j^{-1}$ is  immediate since uniform convergence of
$\{f_n\} \subset C^0(M)$ clearly implies the narrow convergence of $\{f_n \lambda\}$ to $f \lambda.$
\qed

\medskip
We shall now prove that the zeroes of $X$ are the critical points of a certain functional.
Let ${\cal B}_1^+$ be the open subset of ${\cal B}_1$ defined
by $${\cal B}_1^+ = \{f \in {\cal B}_1 \: : \inf_{x \in M} f(x) > 0 \}$$
and let 
$J = J_V : {\cal B}_1^+ \to \RR$ be the $C^{\infty}$ {\em free energy} function defined by
\beq
\label{eq:defJ}
J(f) = \frac{1}{2}  \la Vf, f \ra_{\lambda}  + \la  f, \log(f) \ra_{\lambda} 
\eeq
\brem {\rm It has been pointed to us by Florent Malrieu that the free
  energy $J$ occurs naturally in the analysis of certain non linear
  diffusions used in the modeling of granular flows (see Carillo,
  McCann and Villani (2003), Malrieu (2001)); and by J. Hofbauer that
  a finite dimensional version of $J$ appears in the analysis of some
  ordinary  differential equations in evolutionary game theory (see
  Hofbauer (2000)).}
\erem

The following proposition shows that the zeroes of $X$ are exactly the critical points of $J$ and have the same type (i.e., sinks or saddles). 
\bprop
\label{th:spectrum} Given $f \in {\cal B}_1^+,$ let 
$\mathsf{T}(f) : C^0(M) \to {\cal B}_0$  be the  operator defined
by
\beq \label{eq:defT}
\mathsf{T}(f)h = f h - \la  f, h \ra_{\lambda} f. \eeq
One has
\bdes
\label{eq:D2J}
\iti $\forall u,v \in {\cal B}_0$ $$D^2 J(f)(u,v) =   \la u,v \ra_{1/f} + \la Vu, v \ra_{\lambda} =  \la (Id + \mathsf{T}(f) \circ V)u, v \ra_{1/f}.$$
\itii ${\cal B}_0$ admits a direct sum decomposition 
$${\cal B}_0 = {\cal B}_0^u(f) \oplus {\cal B}_0^c(f) \oplus {\cal B}_0^s(f)$$
where
\bdes
\ita  ${\cal B}_0^u(f), {\cal B}_0^c(f), {\cal B}_0^s(f)$ 
are closed subspaces invariant under $(Id + \mathsf{T}(f) \circ V);$
\itb ${\cal B}_0^c(f) = \{u \in {\cal B}_0: \:  (Id + \mathsf{T}(f) \circ V)u = 0 \}$ and $Id + \mathsf{T}(f) \circ V$ restricted to ${\cal B}_0^{u}(f)$ or ${\cal B}_0^s(f)$ 
is an isomorphism;
\itc  Both ${\cal B}_0^u(f)$ and ${\cal B}_0^c(f)$ have finite dimension;
\itd  The bilinear form $D^2J(f)$ restricted to ${\cal B}_0^u(f)$
(respectively ${\cal B}_0^c(f)$, respectively  ${\cal B}_0^s(f)$) is definite negative (respectively null, respectively
definite positive).
\edes
\itiii We have $$DJ(f) = 0 \quad \Leftrightarrow \quad X(f) = 0,$$ and
in this case, for all $u \in {\cal B}_0$
$$DX(f)u = - (Id + \mathsf{T}(f) \circ V)u.$$
\edes
\eprop 
\prf $(i)$ For all $u \in {\cal B}_0$
\beq
\label{eq:DJ}
DJ(f)u = \la Vf + \log(f) + 1, u \ra_{\lambda} = \la Vf + \log(f), u \ra_{\lambda}.
\eeq
Therefore 
$$
D^2J(f)(u,v) = \la  Vu + \frac{1}{f}u,v  \ra_{\lambda} = \la Vu,v \ra_{\lambda} + \la u,v \ra_{1/f}
$$
which gives the first expression for $D^2J(f).$
Since for all $u,v \in {\cal B}_0$
\beq
\label{eq:Ksym} \la \mathsf{T}(f) V u, v \ra_{1/f} = \la Vu,v \ra_{\lambda} - \la f,Vu \ra_{\lambda}\la \one,v \ra_{\lambda}
= \la Vu,v \ra_{\lambda}
\eeq
we get the second expression for $D^2J(f).$ \\

$(ii)$ Let $K$ denote the operator $\mathsf{T}(f) \circ V$ restricted to $L^2_0(\lambda).$ Then $K$ is compact (by Lemma \ref{th:Vcompact}), and  self-adjoint with respect to the inner product $\la \cdot,\cdot \ra_{1/f}$ (by equation (\ref{eq:Ksym})).
It then follows, from the spectral theory of
compact self-adjoint operators (see Lang, 1993 Chapters XVII and
XVIII) that 
\bdes
\ita $K$ has at most countably many real
eigenvalues;  \itb The set of nonzero eigenvalues is either finite or
can be ordered as $|c_1| > |c_2| > \ldots > 0$ with $\lim_{i \to \infty}
c_i = 0;$   \itc The family $\{{\cal H}_{c}\}$ of eigenspaces, where $c$
ranges over all the eigenvalues (including $0$) forms an orthogonal
decomposition of $L^2_0(\lambda)$;  \itd Each ${\cal H}_{c}$ has finite
dimension provided $c \neq 0$.
\edes
We now set ${\cal B}_0^c(f) = {\cal H}_1$,
${\cal B}_0^u(f) = \oplus {\cal H}_{d}$ where $d$ ranges over
all eigenvalues $ > 1$ and ${\cal B}_0^s(f) = ({\cal B}_0^c(f) \oplus {\cal B}_0^u(f))^{\perp} \cap {\cal B}_0$. \\

$(iii)$ From (\ref{eq:DJ}),  and by density of $\cal{B}_0$ in $L^2_0(\lambda),$ 
$DJ(f) = 0$ if and only if $Vf + \log(f) \in \RR \one.$
Since $f \in {\cal B}_1,$ this is equivalent to $f = \xi(Vf).$
Now, 
\beq
\label{eq:DX}
DX(f) = - Id  - \mathsf{T}(\xi(Vf)) \circ V
\eeq
Hence
$DX(f) = -  Id + \mathsf{T}(f) \circ V$ when $X(f) = 0.$
\qed

\medskip
Let $f \in X^{-1}(0)$, or equivalently $\mu = f \lambda \in
\mathsf{Fix}(\Pi).$ We say that $f$ (respectively, $\mu$) is  a {\em
  nondegenerate} zero or {\em equilibrium} of $X$ (respectively, a
       {\em  nongenerate fixed point} of $\Pi$)  if the space ${\cal
	 B}_0^c(f)$ in the above decomposition reduces to zero.  The
       {\em index} of $f$ (respectively, $\mu$) is defined to be  the
       dimension of ${\cal B}_0^u(f).$

A nondegenerate zero of $X$ (fixed point of $\Pi$)  is called a {\em
  sink} if it has zero index  and  a {\em saddle} otherwise.
\\

Let  $C^k_{sym}(M \times M)$, $k \geq 0$ denote the Banach space of $C^k$ symmetric functions $V : M \times M \to \RR,$ endowed with the topology of $C^k$ convergence. 
The following theorem gives some sense to the hypothesis (made in theorems \ref{th:degree}, \ref{th:main2} and \ref{th:cormain2}  below) that fixed points of $\Pi$ are nondegenerate. However we wont make any other use of  this theorem. The proof is given in the appendix (section \ref{sec:appendix}). 
\bthm
\label{th:generic}
Let ${\cal G}$ denote the  set of $V \in C^k_{sym}(M \times M)$ such that 
$\Pi_V$ has nondegenerate fixed points. Then ${\cal G}$ is open and dense.
\ethm
\brem
{\rm 
The key argument that will be used  in the proof of the genericity Theorem \ref{th:generic} is Smale's infinite-dimensional version of Sard's theorem  for Fredholm maps.   This result by  Smale is also at the origin of the Brouwer degree theory for Fredholm maps initially developed by Elworthy and Tromba  (1970). A  consequence of this degree theory (applied to $X$) is   the following result 
\bthm
\label{th:degree}
Suppose that  every $\mu^* \in \mathsf{Fix}(\Pi)$ is nondegenerate.
Let $C_k$, $k \geq 0$ denote the number of fixed point for $\Pi$ having index $k.$
Then
$$\sum_{k \geq 0} (-1)^k C_k = 1.$$
\ethm
}
\erem
\subsection{Self-repelling diffusions}
A function  $K : M \times M \to \RR$  is called 
a {\em Mercer kernel}, if $K$ is continuous, symmetric and defines a positive operator in the sense that 
$$\la Kf , f \ra_{\lambda} \geq 0$$ for all $f \in  L^2(\lambda).$

If, up to an additive constant\footnote{The dynamics (\ref{eq:sde}) is
 unchanged if one replace $V(x,y)$ by $V(x,y) + \beta.$}, $V$
 (respectively, $-V$) is a Mercer kernel, we call $\{X_t\}$ (given by
 (\ref{eq:sde}))  a {\em self-repelling} (respectively, {\em
 self-attracting} process). 
The following result and the examples below give some  sense to this
 terminology (see in particular examples \ref{ex:sphere},
 \ref{ex:laplace} and \ref{ex:compmon}).
\bthm
\label{th:repul1}
Suppose that, up to an additive constant, $V$ is a Mercer kernel. Then 
\bdes
\iti $J = J_V$ is strictly convex,
\itii 
 ${\mathsf {Fix}}(\Pi)$ reduces to a
singleton $\{\mu^*\}$ and $\lim_{t \to \infty} \mu_t = \mu^*$ almost
surely. If we furthermore assume that hypothesis \ref{hyp:occ} holds,
then $\mu^* = \lambda$.
\edes
\ethm
\prf follows from the definition of $J$,   Proposition
 \ref{th:spectrum} and  Theorem \ref{th:main}. \qed

\bex
\label{ex:mercer}
{\rm Let $C$ be a metric space, $\nu$ a probability over $C$ and $G : M  \times C \to \RR$ a continuous bounded function. Then 
$$K(x,y) = \int_C G(x,u) G(y,u) \nu(du)$$ is a Mercer kernel.
Indeed $K$ is clearly continuous, symmetric,  and  
$$\la Kf, f \ra_{\lambda} = \int_C \left ( \int_M G(x,u)f(x)\lambda(dx)\right )^2 \nu(du) \geq 0.$$ Note that when $C = M$ and $\nu = \lambda$ then
$K = G^2$ as an operator on $L^2(\lambda).$
}
\eex
\bex
\label{ex:sphere} {\bf (i)}
{\rm Let $M = S^d \subset \RR^{d+1}$ be the unit sphere of $\RR^{d+1}$
and let $K(x,y) = \la x, y \ra = \sum_{i = 1}^{d+1} x_i y_i$. Then $K$
is a Mercer kernel  (take $C = \{1, \ldots d+1\},$ $\nu$ the uniform
measure on $C$, and $G(i,x) = \sqrt{d+1} \times x_i$).
}
\eex
\bex 
\label{ex:laplace}{\rm Let $\Delta$ denote the Laplacian on $M$
and $\{K_t(x,y)\}$ the Heat kernel of $e^{\Delta t}$. Fix
$\tau > 0$ and let $K = K_{\tau}.$  The function  $G(x,y) =   K_{{\tau}/2}(x,y)$ is
a  symmetric  $C^{\infty}$ Markov kernel  so that
$K$ is Mercer kernel in view of the  example \ref{ex:mercer} (take $C = M$ and $\nu =
\lambda$). 
}
\eex
\bex 
\label{ex:laplace2}{\rm The example above can be generalized as follows. Let
$\{P_t\}_{t \geq 0}$ be a continuous time Markov semigroup reversible
with respect to some probability measure $\nu$ on $M.$ Assume that
$P_t(x,dy)$ is absolutely continuous with respect to $\nu$ with
smooth density $K_t(x,y)$. Then $K(x,y) = K_{\tau}(x,y)$ is a Mercer kernel.}
\eex
\bex
\label{ex:torus} {\bf (i)}  
{\rm Let $M = T^d = \RR^d/(2\pi \ZZ)^d$  be the flat
$d$-dimensional torus, and let  $\kappa : T^d \to \RR$ be an even
(i.e. $\kappa(x) = \kappa(-x)$) continuous function. Set 
\beq
K(x,y) = \kappa(x-y). \eeq
Given $k \in \ZZ^d,$ let
\beq 
\label{eq:fourco}
\kappa_k = \int_{T^d} \kappa(x) e^{-ik\cdot x} \lambda(dx) \eeq
be the $k$-th Fourier coefficient of $\kappa$. Here
$k\cdot x = \sum_{i = 1}^d k_i x_i$ and $\lambda$ is the normalized
Lebesgue measure on $T^d \sim [0, 2\pi[^d$.
Since $v$ is real and even, $\kappa_{-k} = \kappa_k = \bar{\kappa_k}.$ If we furthermore assume that $$\forall k \in \ZZ^d, \: \kappa_k \geq 0,
$$ then $K$ is a Mercer kernel, since
$$\la Kf, f \ra_{\lambda} = \sum_{k} \kappa_k |f_k|^2$$ for all $f \in L^2(\lambda)$ and  $f_k$ the $k$-th Fourier coefficient of $f.$}
\eex
\bex 
\label{ex:compmon} {\rm  
A  function  $f : [0, \infty[ \to \RR$ is said {\em completely monotonic} if it is  $C^{\infty}$ and, for all $t > 0$ and $k \geq 0,$ 
$$ (-1)^k \frac{d^kf}{dx^k} (t) \geq 0.$$ Examples of such functions are $f(t) = \beta e^{-t/\sigma^2}$
and $f(t) = \beta (\sigma^2 + t)^{-\alpha}$ for $\sigma \neq 0, \alpha, \beta > 0.$
 
Suppose $M \subset \RR^n,$ and 
$K(x,y) = f(||x-y||^2)$ where $f$ is  completely monotonic and $||\cdot||$ is the Euclidean norm on $\RR^n.$
 Then it was proved by 
 Schoenberg (1938)   that $K$ is a Mercer kernel.
}
\eex
\subsubsection*{Weakly self-reppeling diffusions}
When $V$ is not a Mercer kernel but can be written as the difference of two Mercer kernels, it is still possible to give a condition ensuring strict convexity of $J.$

We will need the following consequence  of the so-called Mercer's theorem:
\blem 
\label{th:Mercer}
Let $K$ be  a Mercer kernel. Then 
there exists  continuous symmetric functions 
$G^n : M \times M \to \RR, \, n \geq 1$
such that 
$$K(x,y) = \lim_{n \to \infty} \la G^n_x, G^n_y \ra_{\lambda}$$
uniformly on $M \times M.$ Here $G^n_x$ stands for the function  $u \mapsto G^n(x,u).$
\elem
\prf
The kernel $K$ defines a  compact positive and self adjoint operator on $L^2(\lambda).$  Hence, by the spectral theorem, $K$ has countably (or finitely) many nonnegative eigenvalues  $(c_k^2)_{k \geq 1}$  and  the corresponding eigenfunctions $(e_k)$ can be chosen to form an orthonormal system. 
Furthermore, by Mercer's theorem (see Chap XI-6 in Dieudonn\'e (1972)) 
$K(x,y) = \sum_i c_i^2 e_i(x)e_i(y)$
where the convergence is absolute and uniform. Now
set $G^n_x(y) = G^n(x,y) = \sum_{i = 1}^n c_i e_i(x)e_i(y).$
\qed

\medskip
To a Mercer kernel $K$ we associate the function $D_K : M \times M \to \RR^+$
given by
\beq
\label{defDK} D_K^2(x,y) = \left [ \frac{K(x,x) + K(y,y)}{2} - K(x,y)\right ] \eeq
$$= \lim_{n \to \infty} \frac{1}{2}||G_x^n - G_y^n||^2_{\lambda}$$
where the $(G^n)$ are like in Lemma \ref{th:Mercer}.

Note that $D_K$ is a semi-distance on $M$ (i.e. $D_K$ is nonnegative, symmetric, verifies the triangle inequality, and vanishes on the diagonal).
We let $$\mathsf{diam}_K(M) = \sup_{x,y \in M} D_K(x,y)$$ denote the diameter of $M$ for $D_K.$

Another useful quantity is 
$$K(x,x) =  \lim_{n \to \infty} ||G_x^n||^2_{\lambda}.$$
We let 
$$\mathsf{diag}_K(M) = \sup_{x \in M} K(x,x).$$
\brem
{\rm 
Notice that there is no obvious way to compare $\mathsf{diam}_K(M)$ and $\mathsf{diag}_K(M).$ For instance,
If $K$ is the kernel given in example \ref{ex:compmon},
then $$\mathsf{diam}_K(M) = f(0) - f(\sup_{x,y} ||x-y||^2) \leq \mathsf{diag}_K(M) = f(0),$$  
while, $$\mathsf{diam}_K(M) = 2 > \mathsf{diag}_K(M) = 1$$
with  $K$ the kernel  given in example \ref{ex:sphere}.
}
\erem
\bthm
\label{th:repul2}
Suppose that, up to an additive constant, 
\beq V = V^+ - V^- \eeq
where $V^+$ and $V^-$ are Mercer kernels.

If   $\mathsf{diam}_{V^-}(M)< 1,$ or $\mathsf{diag}_{V^-}(M) < 1,$ then the  conclusions of theorem \ref{th:repul1}
hold. 
\ethm
\prf  First note that 
 $J_V(f) = \frac{1}{2} \la V^+ f, f \ra + J_{-V^-}(f),$ and  since  $f \mapsto \la V^+ f, f \ra_{\lambda}$ is convex,  it suffices to prove that $J_{-V^-}$ is strictly convex.
We can therefore assume, without loss of generality, that $V^+ = 0.$ Or, in other words, that $- V$ is a Mercer kernel. 
We proceed in two steps.

{\em Step 1:} We suppose here that $V(x,y) = - \la G_x, G_y \ra_{\lambda}$
for some  continuous symmetric function $G:(x,u) \mapsto G_x(u).$ 
By Proposition \ref{th:spectrum}, proving that $D^2J_V(f)$ is definite positive reduces to show that   $Id +  \mathsf{T}(f) V =  Id -  \mathsf{T}(f) G^2$ has eigenvalues $ > 0$, or equivalently, that  $\mathsf{T}(f) G^2$ has eigenvalues $< 1.$

Let $\lambda$ be an eigenvalue for  $\mathsf{T}(f) G^2$ and $u \in {\cal B}_0$ a corresponding eigenvector. Set $v = Gu.$ Then 
$$\mathsf{T}(f) Gv = \lambda u.$$ This implies that  $v \neq 0$ (because $u \neq 0$) and  that 
\beq 
\label{eq:tvg}
G \mathsf{T}(f) G v = \lambda v. 
\eeq
 Thus, using the fact that $G$ is symmetric,
$$ \la \mathsf{T}(f) G v, G v \ra_{\lambda} =  \lambda   ||v||^2_{\lambda}.$$
That is
\beq
\label{eq:lamb1}
{\mathsf {Var}}_{f}(Gv) =  \lambda  ||v||^2_{\lambda}.
\eeq
where 
\beqar
{\mathsf {Var}}_{f}(u) & =  &  \la \mathsf{T}(f) u, u \ra_{\lambda}
\label{eq:defvar}\\
&=& \int_M u^2(x) f(x) \lambda(dx) - \left(\int_M u(x)f(x)\right)^2 \lambda(dx) 
\non\eeqar
Now
\beq {\mathsf {Var}}_{f}(Gv) = \frac{1}{2} \int_{M \times M}
(Gv(x) - Gv(y))^2 f(x) f(y) \lambda(dx) \lambda(dy). \eeq
On the other hand,  
\beqarr (Gv(x) - Gv(y))^2 &=& \la  G_x - G_y, v \ra_{\lambda}^2 \\
&\leq& \|G_x - G_y\|^2 \|v\|^2 = 2 (D_{-V}(x,y))^2 \|v\|^2_{\lambda}. \eeqarr
Thus
\beq
\label{eq:lamb2}
 \mathsf{Var}_{f}(Gv) \leq (\mathsf{diam}_{-V})^2  \|v\|^2_{\lambda}.
\eeq
Combining (\ref{eq:lamb1}) and (\ref{eq:lamb2}) leads to $\lambda < (\mathsf{diam}_{-V})^2 < 1.$ 

To obtain the second estimate, observe that (by (\ref{eq:defvar}))
 $${\mathsf {Var}}_{f}(Gv)
\leq \int_M (\la G_x, v \ra)^2 f(x) \lambda (dx)$$
$$ \leq ||v||_{\lambda}^2  \int ||G_x||^2 f(x)\lambda(dx) \leq \mathsf{diag}_{-V}(M) ||v||_{\lambda}^2.$$
{\em Step 2:} In the general case, by lemma \ref{th:Mercer}, we have 
$V(x,y) = \lim_{n \to \infty} V^n(x,y)$ uniformly on $M \times M$
where $V^n(x,y) = - \la G^n_x, G^n_y \ra_{\lambda}.$

Hence, assuming $\mathsf{diam}_{-V}(M) < 1,$ we get that  $\mathsf{diam}_{-V^n}(M) < 1$  for $n$ large enough.
Then, by {\em step 1}, there exists $\alpha > 0$ such that 
 $$D^2J_{V^n}(u,u) = \la u + \mathsf{T}(f) V^n u , u \ra_{1/f} \geq \alpha ||u||^2_{1/f}$$ for all $u \in {\cal B}_0.$
Passing to the limit when $n \to \infty$ leads to  $$D^2J_{V}(u,u) \geq \alpha ||u||_{1/f}^2.$$
The proof of the second estimate is similar.
\qed

\paragraph{Example \ref{ex:sphere} (ii), (continued)} {\rm
Suppose $M = S^d \subset \RR^{d+1}$ and $$V(x,y) = a \times \la x, y \ra
= a \times\sum_{i = 1}^{d+1} x_i y_i$$
for some $a \in \RR.$ The kernel $K = \mathsf{sign}(a) V$ is a Mercer kernel, and  
$\mathsf{diag}_{K}(M) = |a|.$ Hence, by Theorem \ref{th:repul2},
$\mu_t \to \lambda ~a.s$ for $a > -1.$

This condition is far from being sharp since  it actually follows from Theorem 
 4.5  in (BLR) that 
$$a \geq -(d+1) \quad \Leftrightarrow \quad \mu_t \to \lambda  \quad a.s.$$}

\paragraph{Example \ref{ex:torus} (ii), (continued)}
{\rm Let $v$ be an even $C^3$ real valued function defined on the flat d-dimensional torus (see example \ref{ex:torus}) and  $$V(x,y) = v(x-y).$$

As a consequence of theorem \ref{th:repul2} we get the following result which generalizes largely  Theorem 4.14 of (BLR). It also corrects a mistake in the proof of this theorem.
\bprop
Let $(v_k)_{k \in \ZZ^d}$ denote the Fourier coefficients of $v$ as defined by (\ref{eq:fourco}). Assume that
$$\sum_{k \in \ZZ^d \setminus \{0\} } \inf(v_k, 0) > - 1.$$
Then $\mu_t \to \lambda$ almost surely. 
\eprop
\prf
 Integrating by part $3$ times, and using the fact that $v\in
C^3$, proves that for all $k \in \ZZ^d$,  $|v_k| \leq
\frac{C}{\|k\|^3}$, where $||k|| = \sup_i |k_i|$ and $C$ is some
positive constant. Hence the Fourier series
$$v_n(x) = \sum_{\{k \in \ZZ^d \: : ||k|| \leq n\} } {v_k} 
e^{ik\cdot x}$$
congverges uniformly to $v.$

Set $$v^{-}(x) = -  \sum_{\{k \in \ZZ^d \setminus \{0\}  : \: v_k < 0\}} v_k e^{ik\cdot x}.$$ 
Then $v(x) = v^+(x) - v^-(x) + v_0$ where $V^+(x,y) = v^+(x-y)$ and $V^-(x,y) = v^-(x-y)$ are  Mercer kernels.
Clearly, $$\mathsf{diag}_{V^-}(T^d) = v^-(0) = - \sum_{\{k \neq 0 : \: v_k < 0\}} v_k$$ and the result follows from theorem \ref{th:repul2}. 
\qed
\subsection{Self-attracting diffusions}
\label{sec:nonconvex}
The results of this section are motivated by the analysis of
self-attracting diffusions (i.e., $-V$ is a Mercer kernel), but apply
to a more general setting.

\smallskip
Recall that $\mu^*\in {\mathsf{Fix}(\Pi)}$ is a sink if $\mu^*$ is
nondegenerate and has zero index (thus it corresponds to a nondegenerate local
minimum of $J$). We denote by $\mathsf{Sink}(\Pi)$ the set of sinks.

\smallskip
The following result is proved in section \ref{sec:attract}.
\bthm
\label{th:main2}
Let $\mu^* \in {\mathsf{Sink}(\Pi)}$.  Then
$$\P[\lim_{t \to \infty} \mu_t = \mu^*] > 0.$$
\ethm

The next theorem is a converse to Theorem \ref{th:main2}  under a  supplementary condition on $V$ that we now explain.

From the spectral theory of compact self-adjoint operators (see e.g~ Lang, 1993, Chapters XVII and XVIII)
$L^2(\lambda)$ admits an orthogonal decomposition invariant under $V$
$$L^2(\lambda) = E^0_V \oplus E^+_V \oplus E^-_V$$
where $E^0_V$ stands for the kernel of $V$ and $V$ restricted to  $E^+_V$ (respectively, $-V$ restricted to  $E^-_V$)
is a positive operator. 

Let $\pi_+$ and $\pi_-$ be respectively  the orthogonal projections from
$L^2(\lambda)$ onto  $E^+_V$ and  $E^-_V.$ Set
\beq
\label{defV+}
V_+ = V \circ \pi_+ \quad\mbox{ and } \quad
V_- = - V \circ \pi_-.
\eeq
So that $V = V_+ - V_-.$
\begin{hypothesis} [Occasional assumption 2]
\label{hyp:occ2}
 $V_+$ and $V_-$ are Mercer kernels. 
\end{hypothesis}

Recall that $\mu^*\in {\mathsf{Fix}(\Pi)}$ is a saddle if $\mu^*$ is
nondegenerate and has positive index. The following theorem is proved
in section \ref{sec:nonconv}.

\bthm
\label{th:main3}
Assume that hypothesis \ref{hyp:occ2} holds.
Let $\mu^* \in {\mathsf{Fix}(\Pi)}$ be  a saddle.  Then
$$\P[\lim_{t \to \infty} \mu_t = \mu^*] =  0.$$
\ethm

\bcor 
\label{th:cormain2}
Suppose  that hypothesis \ref{hyp:occ2} holds and that every $\mu^* \in {\mathsf{Fix}(\Pi)}$ is nondegenerate. 
Then there exists a random variable $\mu_\infty$ such that
\bdes \iti $\lim_{t\to\infty}\mu_t = \mu_\infty$ a.s.
\itii $\P[\mu_\infty\in\mathsf{Sink}(\Pi)]=1$ and
\itiii For all $\mu^*\in\mathsf{Sink}(\Pi)$, 
$$\P[\mu_\infty = \mu^*] >0 . $$
\edes\ecor
\prf follows from Theorems \ref{th:main}, \ref{th:main2} and \ref{th:main3}. \qed
\subsection{Localisation}
In this section, we assume that  hypothesis \ref{hyp:occ} holds. In this
case, $\lambda$ is always a fixed point for $\Pi$, hence a possible limit point for $\{\mu_t\}.$  We will say that
the self-interacting diffusion ``localizes'' provided $\P[\mu_t\to\lambda]=0$.
We have already seen (see Theorems \ref{th:repul1} and \ref{th:repul2}) that
self-reppelling diffusions and weakly self-attracting diffusions
never localize.
\bthm
\label{th:loc1}
Suppose that hypothesis \ref{hyp:occ} holds. Let
\beq\rho(V) = \inf \{ \la  Vu, u \ra_{\lambda}: \: u \in L^2_0(\lambda), \, \|u\|_{\lambda} = 1 \}.\eeq
Assume that $\rho(V)>-1$, then
\beq \P[\lim_{t \to \infty} \mu_t = \lambda] > 0 . \eeq
Assume that $\rho(V)<-1$ and that hypothesis \ref{hyp:occ2} holds, then 
\beq \P[\lim_{t \to \infty} \mu_t = \lambda] = 0 . \eeq
\ethm
\prf  Under hypothesis \ref{hyp:occ},  $\xi(V\lambda) = 1.$ Then, by Proposition \ref{th:spectrum}  
$$D^2J(1)(u,v) = - \la DX(1)u,v \ra_{\lambda} = \la u + Vu, v \ra_{\lambda}.$$
The result then follows from Theorems \ref{th:main2} and \ref{th:main3}. \qed

\paragraph{Example \ref{ex:torus} (iii), (continued).}
With $V$ as in example \ref{ex:torus} (ii),
$$\rho(V) = \inf_{k \in  \ZZ^d \setminus \{0\}} v_k.$$

\paragraph{Example \ref{ex:laplace} (ii), (continued).}
Suppose  $V(x,y) =  a K_{\tau}(x,y)$ for some  $a \leq 0$ and $\tau > 0,$ where $\{K_t\}_{t > 0}$ is the Heat kernel of $e^{\Delta t}.$ 
Then $\rho(V) = ae^{-\lambda \tau}$ where
 $\lambda$ is the smallest non zero eigenvalue of $\Delta.$ Note that
there exist numerous estimates of $\lambda$ in terms of the geometry
of $M$.
\section{Review of former results}
\label{sec:background}
We recall here some notation  and results from (BLR) on which rely our
analysis. There is no assumption in this section that $V$ satisfies
one of the hypotheses \ref{hyp:main} or \ref{hyp:occ}. The only
required assumption is that $V$ is smooth enough, say\footnote{see
(BLR) for a more precise assumption} $C^3$.

The map $\Pi$ defined by  (\ref{eq:PI}) extends to a map  $\Pi :{\cal M}(M) \to {\cal P}(M)$ given by the same formulae. 
Let $F : {\cal M}_s(M) \to {\cal M}_s(M)$ be the vector field defined
by 
\beq \label{eq:ode}
F(\mu) =  - \mu + \Pi(\mu), \eeq
Then (see (BLR), Lemma 3.2)  $F$ induces a $C^{\infty}$ flow
$\{\Phi_t\}_{t \in \RR}$ on ${\cal M}_s(M).$

The  {\em limiting  dynamical system associated to $V$} is the mapping 
\beq \begin{array}{lllll}
\Psi &:&  \RR \times {\cal P}_w(M)  &\to& {\cal M}_w(M), \\
&& (t,\mu)  &\mapsto& \Psi_t(\mu)=\Phi_t(\mu). \end{array}\eeq
Because $\Phi$ is a flow, $\Psi$  satisfies the flow property:
\beq \Psi_{t+s}(\mu) = \Psi_t \circ \Psi_s(\mu) \eeq
for all $t,s \in \RR$ and $\mu \in {\cal P}(M) \cap \Phi_{-s}({\cal
P}(M))$. Furthermore, (see Lemmas 3.2 and 3.3 of (BLR)) $\Psi$ is
continuous and leaves  ${\cal P} (M)$ positively invariant:
\beq \Psi_t({\cal P} (M)) \subset {\cal P}(M) \hbox{ 
for all }t \geq 0. \eeq
The key tool for analyzing self-interacting diffusion is  
Theorem \ref{th:limitset} below (Theorem 3.8 of (BLR)), according to
which, the long term behavior of the sequence $\{\mu_t\}$ can be
described in term of certain invariant sets for $\Psi$.
Before stating this theorem, we first recall some definitions from
dynamical systems theory.
\subsubsection*{Attractor free sets and the Limit set theorem}
A subset $A \subset {\cal P}_w(M)$ is said to be {\em invariant} for
$\Psi$ if $\Psi_t(A) \subset A$ for all $t \in \RR$. Let  $A$ be an
invariant  set for $\Psi$. Then $\Psi$ induces a flow  on $A$
$,\Psi|A$  defined by taking the restriction  of $\Psi$ to $A$. That
is $(\Psi|A)_t = \Psi_t|A.$

Given an invariant  set $A$, a set $K \subset A$ is called an {\em
attractor} (in the sense of Conley (1978)) for $\Psi|A,$ if it is
compact, invariant and has a neighborhood $W$ in $A$ such that
\beq \lim_{t \to \infty} \dist_w(\Psi_t(\mu), K) = 0 \eeq
uniformly in $\mu \in W$. Here $\dist_w$ is any metric on
$\mathcal{P}_w(M)$.

An attractor $K\subset A$ for $\Psi|A$ which is different from
$\emptyset$ and $A$ is called {\em proper}.
An {\em attractor free set} for $\Psi$ is a nonempty compact invariant
set $A \subset {\cal P}_w(M)$ with the property that $\Psi|A$ has no
proper attractor. Equivalently, $A$ is a nonempty compact connected
invariant set such that $\Psi|A$ is a {\em chain-recurrent} flow
(Conley, 1978).

\brem 
\label{gendyn}
{\rm
The definitions (invariant sets, attractors, attractor free sets)
given here for $\Psi$ extend obviously to any  (local) flow on a
metric space. This will be used below.}
\erem

The {\em limit set} of $\{\mu_t\}$ denoted $L(\{\mu_t\})$ is the set
of limits of convergent sequences $\{\mu_{t_k}\}$, $t_k \to
\infty$. That is 
\beq L(\{\mu_t\}) = \bigcap_{t \geq 0} \overline {\{\mu_s \,: s \geq
t\}}\eeq
where $\bar{A}$ stands for the closure of $A$ in ${\cal P}_{w}(M).$

\bthm [(BLR), Theorem 3.8]
\label{th:limitset}
With probability one $L(\{\mu_t\})$
is an attractor free set of $\Psi.$
\ethm

This result allows,  in various situations, to characterize exactly
the asymptotic of $\{\mu_t\}$ in term of the potential $V$ and the
geometry of $M$. We refer the reader to (BLR) for several  examples
and further results. Amongst the general consequences of Theorem
\ref{th:limitset} the two following corollaries will be useful here.
\bcor
\label{th:limitsetcor}
Let $A \subset {\cal P}_w(M)$ be an attractor and
\beq B(A) = \{\mu \in {\cal P}_w(M) \: : \lim_{t \to \infty}
\dist_w(\Psi_t(\mu), A) = 0\} \eeq
its basin of attraction. Then the events
\beq \{ L(\{\mu_t\}) \bigcap B(A)  \neq \emptyset\}\quad \mbox{ and } 
\quad \{L(\{\mu_t\}) \subset A\} \eeq
coincide almost surely. \ecor

For a proof see ((BLR), Proposition 3.9).
\bcor \label{th:smoothdens}
With probability one, every point $\mu^* \in L(\{\mu_t\})$ can be
written as
\beq \mu^* = \int_{{\cal P}_w(M)} \Pi(\mu) \rho(d\mu) \eeq
where $\rho$ is a Borel  probability measure over  ${\cal P}_{w}(M).$
In particular, if $V$ is $C^k$ then $\mu^*$ has a $C^k$ density with
respect to $\lambda$. \ecor

This last result follows from Corollary \ref{th:limitsetcor} as
follows: Let
\beq C_{\Pi}({\cal P}_w(M)) = \left\{ \int_{ {\cal P}(M)} \Pi(\mu)
\rho(d\mu): \: \rho \in {\cal P}({\cal P}_w(M))\right\} \eeq
where ${\cal P}({\cal P}_w(M))$ is the set of Borel probability
measures over ${\cal P}_w(M).$ It is not hard to prove that
$C_{\Pi}({\cal P}_w(M))$ contains a global attractor for $\Psi;$ that
is  an attractor whose basin is ${\cal P}_w(M)$. Hence $ L(\{\mu_t\})
\subset C_{\Pi}({\cal P}_w(M))$ by Corollary \ref{th:limitsetcor}. For
details see ((BLR), Theorem 4.1).

\section{Convergence of $\{\mu_t\}$ toward $\mathsf{Fix}(\Pi)$}
\label{sec:proof1}
This section is devoted to the proof of Theorem \ref{th:main}.  Hypothesis \ref{hyp:main} is implicitly assumed. 
\subsection{The flow induced by $X$}
Recall that ${\cal B}_1^+ = \{f \in {\cal B}_1:  \: f > 0 \},$ where
$\cB_1=\{f\in C^0(M):\:\int f~d\lambda = 1 \}$.
\bprop
\label{th:Jlyapou}
The vector field $X$ given by (\ref{defX})
induces a global smooth flow $\Phi^X = \{\Phi_t^X\}$ on ${\cal B}_1.$ Furthermore,
\bdes
\iti $\Phi_t^X(f) \in {\cal B}_1^+$ for all $t \geq 0$ and $f \in {\cal B}_1^+.$ 
\itii For all $f \in {\cal B}_1^+$ and  $t > 0,$  $J(\Phi_t^X(f)) < J(f)$ if $f$ is not an equilibrium. 
\edes 
\eprop
\prf
The vector field $X$ being smooth, it induces a smooth local flow $\Phi^X$ on ${\cal B}_1.$ To check that this flow is global observe that
$$||- f + \xi(Vf)||_{L^1(\lambda)} \leq ||f||_{L^1(\lambda)} + 1.$$
Hence, by standard results,  the differential equation 
$$\frac{df}{dt} = - f + \xi(Vf)$$ generates a smooth global flow on $L^1(\lambda)$
whose restriction to ${\cal B}_1$ is exactly $\Phi.$ 

$(i)$ For $f \in {\cal B}_1^+$, $||Vf||_{\infty} \leq ||V||_{\infty}.$ Thus
$X(f)(x) \geq -f(x) + \delta$ for all $x \in M,$  where $\delta =  e^{-2||V||_{\infty}}.$
It follows that $\Phi_t^X(f)(x) \geq e^{-t}(f(x) - \delta) + \delta \geq \delta(1-e^{-t}) > 0$ for all $t >  0.$

$(ii)$ For $f \in {\cal B}_1^+,$ let $K_f : {\cal B}_1^+ \rar \RR$
be the ``free energy'' function associated to the potential $Vf$
$$K_f(g) =  \la Vf,g \ra_{\lambda} + \la g, log(g) \ra_{\lambda}.$$
 The function $K_f$ is a  $C^{\infty},$ strictly convex function  and  reaches its global minimum at the ``Gibbs'' measure $\xi(Vf).$
Indeed, a direct computation shows that for $h\in\cB_0$,
$$DK_f(g).h  = \la log(g) + Vf, h \ra_{\lambda}$$
and for $h$ and $k$ in $\cB_0$,
$$D^2K_f(g)(h,k) = \la h,k \ra_{1/g}.$$
Thus $DK_f(g) = 0$ if and only if $g = \xi(Vf)$
and
$D^2K_f(g)$ is positive definite for all $g.$
Then, since
\beq
DK_f(g).[g - \xi(Vf)] = [DK_{f}(g) - DK_{f}(\xi(Vf))].[g - \xi(Vf)],
\eeq
by strict convexity we then deduce that
\beq\label{eq:convex}
DK_f(g).[g - \xi(Vf)] \geq 0,\eeq
with equality if and only if $g = \xi(Vf).$

Now observe that $DJ(f) = DK_f(f).$
Hence, by (\ref{eq:convex}) $$DJ(f).X(f) \leq 0$$ with equality if and
only if $X(f) = 0$. This proves (ii). \qed
\subsection{Proof of Theorem \ref{th:main}}
\blem
\label{th:icont}
The map $i:  C_{\Pi}({\cal P}_w(M))  \to {\cal B}_1^+ \subset C^0(M)$ defined by $i(\mu) = \frac{d\mu}{d\lambda}$
is continuous.
\elem
\prf 
Let $\mu_n = \int_{{\cal P}(M)} \Pi(\nu)\rho_n(d\nu) \in  C_{\Pi}({\cal P}_w(M))$ be such that $\mu_n \rar \mu$ (for the narrow topology).
By Lipschitz continuity of $V$, the family $\{\xi(V\nu), \, : \, \nu \in {\cal P}(M)\}$ is uniformly bounded and equicontinuous. Hence the sequence of densities
$f_n = \int_{{\cal P}(M)} \xi(V\nu)\rho_n(d\nu), n \geq 0$ is  uniformly bounded and equicontinuous. By Ascoli theorem it is relatively compact in $C^0(M).$ 
It easily follows that $f_n \rar f = \frac{d\mu}{d\lambda}$ in
$C^0(M).$ \qed

\blem
\label{th:conjug}
 Let $K \subset {\cal P}_w(M)$ be a compact invariant set for $\Psi.$
Then for all $\mu \in K$ and $t\in\RR$,
$$\Phi_t^X \circ i(\mu) =  i \circ \Psi_t(\mu).$$
\elem
\prf Note that for all $\mu \in C_{\Pi}(\cP(M))$, $X \circ i(\mu) =  i
\circ F(\mu)$ from which the result follows since $K\subset
C_\Pi(\cP(M))$ is invariant. \qed

\medskip
To shorten notation, we set here $L = L(\{\mu_t\}).$ Recall that $L \subset C_{\Pi}(\cP(M))$ (Corollary \ref{th:smoothdens}) and that $L$ is attractor free for $\Psi$ (Theorem \ref{th:limitset}). 
\blem \label{th:conjug2}
$i(L)$ is an attractor free set for $\Phi$.
\elem
\prf
This easily  follows from the continuity of $i$ (Lemma
\ref{th:icont}), compactness of $L$ and the  conjugacy property (Lemma
\ref{th:conjug}) (compare to Corollary 3.10 in (BLR)).  \qed
\bcor  \label{th:sard1}
$i(L)$ is a connected subset of $X^{-1}(0).$
\ecor

Before proving this corollary, remark that it  implies Theorem \ref{th:main} since 
$i^{-1}(X^{-1}(0)) = \mathsf{Fix}(\Pi).$
\paragraph{Proof of Corollary \ref{th:sard1}:}
The proof of this corollary relies on the following result (Benaim (1999), Proposition 6.4):
\bprop
\label{th:ben}
Let $\Lambda$ be a compact invariant set for a flow $\Theta =
\{\Theta_t\}_{t \in \RR}$ on a metric space $E$. Assume there exists a
continuous function ${\cal V} : E \to \RR$ such that
\bdes
\ita ${\cal V}(\Theta_t(x)) < {\cal V}(x)$ for $x \in E \setminus \Lambda$ and $t > 0$.
\itb ${\cal V}(\Theta_t(x)) = {\cal V}(x)$ for $x \in \Lambda$ and $t \in \RR$.
\edes
Such a ${\cal V}$ is called a {\em Lyapounov function} for $(\Lambda, \Theta).$
If ${\cal V}(\Lambda)$ has empty interior, then every attractor free set $K$ for 
$\Theta$ is contained in $\Lambda.$ Furthermore ${\cal V}|K$ (${\cal V}$ restricted
to $K$) is constant.
\eprop

Set $E = i(L), \, \Theta = \Phi^X|i(L),$  $\Lambda = X^{-1}(0) \cap
i(L)$ and ${\cal V} = J|i(L).$  Then  $\Lambda$ is a compact set
(lemma \ref{th:fixzero}), and ${\cal V}$ is a Lyapounov function for
$(\Lambda, \Theta)$ by Proposition \ref{th:Jlyapou}. By Lemma
\ref{th:conjug2}, $i(L)$ is an attractor free set. Therefore,
to apply Proposition \ref{th:ben}, it
suffices  to check  that $J(X^{-1}(0))$ has empty interior. This is a
consequence of the infinite dimensional version of Sard's theorem for
$C^{\infty}$ functionals proved by Tromba (see Theorem 1 and Remark 7
of Tromba, 1977). Thus Proposition \ref{th:ben} proves that
$i(L)\subset X^{-1}(0)$.

\bthm \label{th:tromba} (Tromba, 1977).
Let ${\cal B}$ be a $C^{\infty}$ Banach manifold, $X$ a $C^{\infty}$
vector field on  ${\cal B}$ and $J : {\cal B} \to \RR$ a $C^{\infty}$
function. Assume that
\bdes
\ita $DJ(f) = 0$ if and only if $X(f) = 0,$
\itb $X^{-1}(0)$ is compact,
\itc For each $f \in X^{-1}(0)$, $DX(f) : T_f {\cal B} \to
T_f {\cal B}$ is a Fredholm operator.
\edes
Then $J(X^{-1}(0))$ has empty interior.
\ethm

The verification that Tromba's theorem applies to the present setting is immediate. Indeed, 
assertion $(a)$ follows from Proposition \ref{th:spectrum} and  assertion $(b)$  from  Lemma \ref{th:fixzero}. Recall that a bounded operaror $T$ from one Banach space $E_1$ to a Banach space $E_2$ is Fredholm if its kernel 
$\mathsf{Ker}(T)$ has finite dimension and its range  $\mathsf{Im}(T)$  has finite codimension. Hence assertion $(c)$ follows from Proposition \ref{th:spectrum}. This concludes the proof of Corollary
\ref{th:sard1}. \qed
\section{Convergence toward sinks}
\label{sec:attract}
The purpose  of this  section is to prove Theorem \ref{th:main2}.
\subsection{The vector field $Y = Y_V$}
\label{sec:defY}
In order to prove theorem \ref{th:main2}, it is convenient to 
introduce a new  vector field 
\beq\begin{array}{lllll}
Y = Y_V &:& C^0(M) &\to& C^0(M)\\
         &&f &\mapsto& -f + V \xi(f) \end{array} \eeq
as well as the stochastic  process $\{V_t\}_{t \geq 0}$
defined by 
\beq V_t = V\mu_{e^t}.\eeq
The reason for this is, roughly speaking, the following.
The measure $\mu_t$ is  singular with respect to $\lambda$,
while $\Phi^X$ is defined on a space of  continuous densities.  This  is  not a  problem if we are  dealing with
{\bf qualitative} properties of $L(\{\mu_t\})$ (like in Theorem \ref{th:main}) since we know (by Corollary \ref{th:smoothdens}) that $L(\{\mu_t\})$ consists of measures having smooth densities. 

Proving Theorem \ref{th:main2}  requires {\bf quantitative}  estimates on
the way  $\{\mu_t\}$ approaches its limit set. We shall do this by  showing  
that $\{V_{t+s}\}_{s \geq 0}$ ``shadows'' at a certain rate the deterministic solution  to the Cauchy problem
$$\dot{f} = Y(f)$$ with initial condition $f_0 = V_t.$  
\blem
\label{th:Ysink}
The vector field $Y$ induces a global smooth flow $\Phi^Y = \{\Phi_t^Y\}$
on $C^0(M).$ Furthermore
\bdes
\iti $V \Phi_t^X(f) = \Phi_t^Y(Vf)$ for all $f \in {\cal B}_1$ and $t \in \RR.$
\itii $V$ maps homeomorphically $X^{-1}(0)$ to $Y^{-1}(0),$  sinks to sinks and saddles to saddles.
\edes
\elem
\prf
The vector field $Y$ is $C^{\infty}$ and sublinear because $||Y(f)||_{\infty} \leq ||f||_{\infty} + ||V||_{\infty}.$ It then  induces a global smooth flow.

$(i)$ follows from the conjugacy $V \circ X = Y \circ V.$ 

$(ii)$. It is easy to verify that $V$ induces an homeomorphism from $X^{-1}(0)$ to $Y^{-1}(0)$ whose inverse is $\xi.$ Let $f \in X^{-1}(0)$  and $g = Vf.$ Then with the notation of proposition \ref{th:spectrum}, 
$DX(f) = - (Id +  \mathsf{T}(f) \circ V)$ 
and $DY(g) = - (Id + V \circ  \mathsf{T}(\xi(g)) = - (Id + V \circ \mathsf{T}(f)).$

For all $\alpha\in\RR$, let 
\beqarr
E^\alpha &=& \{u\in L^2(\lambda), \mathsf{T}(f)Vu=\alpha u\} \\
H^\alpha &=& \{u\in L^2(\lambda), V\mathsf{T}(f)u=\alpha u\}.
\eeqarr
The operators $\mathsf{T}(f)V$ and $V\mathsf{T}(f)$ are compact
operators acting on $L^2(\lambda)$. The adjoint of
$\mathsf{T}(f)V$ is $V\mathsf{T}(f)$. This implies that for
$\alpha\neq 0$, $E^\alpha$ and $H^\alpha$ are isomorphic, with
$V:E^\alpha\to H^\alpha$ having for inverse function
$\frac{1}{\alpha}\mathsf{T}(f)$.
Therefore, if  $f$ is nondegenerate (respectively a sink, respectively a saddle) for $X,$ then  $Vf$ is   nondegenerate (respectively a sink, respectively a saddle) for $Y.$
\qed
\subsection{Proof of Theorem \ref{th:main2}}
 We now follow the line of the proof of Theorem 4.12 (b) in (BLR).
We let  ${\cal F}_t$ denote the sigma field generated by the random variables
$(B^i_s: \: s \leq e^t, \, i = 1 \ldots N).$
\blem \label{th:pta} 
There exists a constant $K$ (depending on $V$) such
that for all $T > 0$ and $\delta>0$,
\beq \P\left[ \sup_{0 \leq s \leq T} \|V_{t+s} - \Phi_s^Y(V_t)\|_{\infty}
\geq \delta | {\cal F}_t\right] \leq \frac{K}{\delta^{d+2}} e^{-t}.\eeq
\elem
\prf Given $t \geq 0$ and $s \geq 0$ let $\vareps_t(s) \in {\cal M}(M)$
be the measure defined by
\beq \vareps_t(s) = \int_t^{t+s} (\delta_{X_{e^r}} - \Pi(\mu_{e^r}))
dr.\eeq
Let us first show
\blem \label{lem:pta}
There exists a constant $K$ (depending on $V$) such
that for all $T > 0$ and $\delta>0$,
\beq \P\left[ \sup_{0 \leq s \leq T} \|V\vareps_t(s)\|
\geq \delta | {\cal F}_t\right] \leq \frac{K}{\delta^{d+2}} e^{-t}.\eeq
\elem
\prf  According to  Theorem 3.6 (i) (a) in (BLR)  there exists
a constant $K$ such that for all $\delta>0$ and $f\in C^\infty(M)$,
\beq \label{3.6} \P\left[ \sup_{0 \leq s \leq T} |\vareps_t(s)f|
\geq \delta | {\cal F}_t\right] \leq \frac{K}{\delta^2}\|f\|_\infty^2
e^{-t}.\eeq
Note that this also holds for all $f\in C^0(M)$ (for a larger constant
$K$) since $f$ can be uniformly approximated by smooth functions.
By compactness of $M$ and Lipschitz continuity of $V,$ there exists
a finite set $\{x_1, \ldots, x_m\} \in M$ such that for all $x \in M$ 
$$|V(x,y) - V(x_i,y)| \leq \frac{\delta}{4T}$$
for some $i \in \{1, \ldots, m\}.$
Therefore 
\beqarr
\sup_{0 \leq s \leq T} \|V\vareps_t(s)\|_{\infty} 
&\leq& \sup_{i = 1, \ldots, m} \sup_{0 \leq s \leq T} \|V\vareps_t(s)(x_i)\| + \delta/2. 
\eeqarr
Hence, 
\beqarr
\P\left[ \left.\sup_{0 \leq s \leq T} \|V\vareps_t(s)\|_{\infty} \geq
\delta\right|\mathcal{F}_t\right]\hskip-100pt&&\\
&\leq& \P\left[\left.  \sup_{i = 1, \ldots, m}  \sup_{0 \leq s \leq T}
|\vareps_t(s)V_{x_i}| \geq \delta /2\right|\mathcal{F}_t\right]\\
&\leq& \frac{4mK \|V\|^2_{\infty}}{\delta^2} \times e^{-t}.
\eeqarr 
Since $M$ has dimension $d,$ $m$ can be chosen to be $m = O(\delta^{-d})$
and the result follows. \qed
\medskip

Note that for all $u \in M$
\beqarr
\frac{dV_t(u)}{dt} &=& -V_t(u) + V(u,X_{e^t}) \\
&=& [V F(\mu_{e^t}) + V(\delta_{X_{e^t}} -
\Pi(\mu_{e^t}))](u). \eeqarr
Thus, using the fact that $V F(\mu) = Y (V\mu)$ we obtain
\beqarr
V_{t+s}(u) - V_t(u) &=& \int_t^{t+s} V F(\mu_{e^r})(u)dr  +
V\vareps_t(s)(u) \\
&=& \int_t^{t+s} Y(V_r)(u)dr  + V\vareps_{t}(s)(u) \\
&=& \int_0^{s} Y(V_{t+r})(u)dr  + V\vareps_t(s)(u) \eeqarr
for all $u \in M.$ In short, 
\beq V_{t+s} - V_t = \int_0^{s} Y(V_{t+r})dr  + V\vareps_t(s).\eeq

Let $v(s) = \|V_{t+s} - \Phi^Y_s(V_t)\|_{\infty}$.
Then for $0 \leq s \leq T$
\beq \label{eq:Gron} v(s) \leq \int_0^s \|Y(V_{t + r}) -
Y(\Phi^Y_r(V_t))\|_{\infty} dr
+ \sup_{0 \leq s \leq T} \| V\vareps_t(s)\|_{\infty}. \eeq
Now, for  $t$, $r \geq 0$ both $V_{t+r}$ and $\Phi^Y_r(V_t)$ lie
in $V {\cal P}_w(M)$ which is a compact subset of $C^0(M)$ (by Lemma \ref{th:Vcompact}). Therefore, by Gronwall's lemma
\beq \label{eq:gronwall} 
\sup_{0 \leq s \leq T} v(s) \leq
e^{LT}  \sup_{0 \leq s \leq T} \| V\vareps_t(s)\|_{\infty} \eeq
where $L$ is the Lipschitz constant 
of $Y$ restricted to  $V {\cal P}_w(M)$.

Then, with the estimate (\ref{eq:gronwall}), Lemma \ref{th:pta} follows from Lemma \ref{lem:pta}. \qed

\medskip
The following lemma  is  Theorem 3.7 of (Benaim, 1999) (see also
Proposition 4.13 of (BLR)) restated in the present context. 
\blem \label{th:blr413}
Let $A \subset C^0(M)$ be an attractor for $\Phi^Y$
with basin of attraction $B(A).$ Let $U \subset B(A)$ be an open set
with  closure $\bar{U} \subset B(A).$ Then there exist positive
numbers $\delta$ and $T$ (depending on $U$ and
$\{\Phi^Y\}$) such that 
\beq \P \left[\lim_{t \to \infty} \mathsf{dist}(V_t, A) = 0\right] \geq 
\left(1 - \frac{K}{\delta^{d+2}} e^{-t}\right)\times \P[\exists s \geq t \: : \: V_s \in U] \eeq
where $K$ is given by Lemma \ref{th:pta} and $\mathsf{dist}(\cdot,\cdot)$
is the distance associated to $\|\cdot\|_{\infty}$. \elem

\blem
\label{th:girsan} Let $\mu \in {\cal P}(M)$, $f = V\mu$ and $U$ a neighborhood of
$f$ in $C^0(M)$. Then for all $t > 0$
\beq \P[V_t \in U ] > 0. \eeq
\elem
\prf Let $\Omega_M$ (respectively, $\Omega_{\RR^N}$) denote the space
of continous paths from $\RR^+$ to $M,$ (respectively, $\RR^N$)
equipped with the topology of uniform convergence on compact
intervals and the associated Borel $\sigma$-field.

Let $B_t = (B^1_t, \ldots, B^N_t)$ be a standard Brownian motion on
$\RR^N.$ We let $\P$ denote the law of $(B_t \: : t \geq 0) \in
\Omega_{\RR^N}$ and $\E$ the associated expectation.

Let $\{W^x_t\}$ be the solution to the SDE
\beq \label{eq:defWx}
dW^x_t = \sum_{i = 1}^N F_i(W^x_t) \circ dB^i_t \: :  W^x_0 = X_0 = x \in M
\eeq
Then $W^x \in \Omega$ is a Brownian motion on $M$ starting at $x.$
Let 
\beq \begin{array}{lll}
M(t) &=& {\displaystyle \exp \left( \int_0^t \sum_i \la  \nabla
V_{\mu_s(W)}(W_s), F_i(W_s) \ra dB^i_s \right.}\\
&& {\displaystyle \left. \hskip40pt - \frac{1}{2} \int_0^t \|
\nabla V_{\mu_s(W)}(W_s)\|^2 ds\right)} \end{array}\eeq
where for all path $\omega  \in \Omega$
\beq \mu_t(\omega) = \frac{1}{t} \int_0^t \delta_{\omega_s} ds. \eeq
Then, $\{M_t\}$ is a martingale with respect to 
$(\Omega_{\RR^N}, \{\sigma(B_s, s \leq t)\}_{t \geq 0}, \P)$ and; by
the transformation of drift formula (Girsanov's theorem) (see section
IV 4.1 and Theorem IV 4.2 of Ikeda and Watanabe (1984))
\beq \P[V_t \in U] = \P[V\mu_{e^t} \in U] = \E[M(e^t)
1_{\{V\mu_{e^t}(W) \in U\}}].\eeq
By continuity of the maps $V : {\cal P}_w(M) \to C^0(M)$ (lemma
\ref{th:Vcompact}) and $\omega \in \Omega_M \mapsto \mu_t(\omega) \in {\cal
P}_{w}(M)$ the set ${\cal U} = \{\omega \in \Omega \: :
V\mu_{e^t}(\omega) \in U\}$ is an open subset of $\Omega_M.$ Its
Wiener measure $\P[W \in {\cal U}] = \P[V\mu_{e^t}(W) \in U]$ is then
positive. This implies that  $\E[M(e^t) 1_{\{V\mu_{e^t}(W) \in U\}}] >
0.$ \qed

\medskip
The proof of Theorem \ref{th:main2} is now clear. Let $\mu^*$ be a
sink for $\Pi.$ Then $V^* = V\mu^*$ is a sink for $Y$ according to Lemma \ref{th:Ysink}, and 
Lemmas \ref{th:blr413} and \ref{th:girsan} imply that
$$\P[ V_t \to V^*] > 0.$$
On the event $\{V_t\to V^*\}$, 
$$L(\{\mu_t\}) \subset \{\mu\in\mathsf{Fix}(\Pi):\:V\mu=V^*\}.$$
Note that $\mu\in\mathsf{Fix}(\Pi)$ with $V\mu=V^*$ implies that
  $\mu=\mu^*$. Therefore, on the event $\{V_t\to V^*\}$, we have
  $\lim_{t\to\infty}\mu_t=\mu^*$. This proves Theorem \ref{th:main2}.

\section{Non convergence towards unstable equilibria}
\label{sec:nonconv}
The purpose of this section is  to prove Theorem \ref{th:main3}. That is 
\beq \label{eqinst}\P[\mu_t\to\mu^*]=0\eeq
provided $\mu^*\in \mathsf{Fix}(\Pi)$ is a nondegenerate unstable
equilibrium and  hypothesis \ref{hyp:occ2} holds.

The proof of this result is somewhat long and technical. For the reader's convenenience we first briefly explain our strategy.

$\bullet$ Set $h_t = V\mu_t.$ To prove that $\mu_t \not \to \mu^*$
we will prove that $h_t \not \to h^*.$
We see $h_t$ as a random perturbation of a deterministic dynamical system induced by a vector field
$\tilde{Y}.$  The vector field $\tilde{Y}$ is introduced in subsection \ref{sec:tidy}. It is defined like  the vector field $Y$ (see section \ref{sec:attract}) but on a  subset 
$\cH^K$ of $C^0(M)$ equipped with a convenient Hilbert space structure (subsection \ref{sec:mercer}).

$\bullet$ The fact that $\mu^*$ is a saddle makes $h^*$ a saddle for $\tilde{Y}.$ According to the stable manifold theorem, the set of points whose forward trajectory (under $\tilde{Y}$) remains close to $h^*$ is a smooth submanifold $W^s_{loc}(h^*)$ of nonzero finite codimension. We construct in subsection \ref{sec:stabman} a ``Lyapounov function'' $\eta$  which increases strictly along forward trajectory of $\tilde{Y}$  off $W^s_{loc}(h^*)$ and vanishes on $W^s_{loc}(h^*).$  

$\bullet$ The strategy of the proof now consists to show that
$\eta(h_t) \not \to 0$ (since $\mu_t \to \mu^*$ implies $\eta(h_t) \to 0.$)
Using stochastic calculus (in $\cH^K$) we derive the stochastic evolution of 
$\eta(h_t)$ (subsection \ref{sec:itocal}) and then prove the theorem
in subsections \ref{sec:lemmafirst} and \ref{sec:lemmasec}.

\medskip
In the different (but related) context of urn processes and stochastic approximations, the idea of using the stable manifold theorem to prove the nonconvergence toward unstable equilibria is due to Pemantle (1990). 
Pemantle's probabilistic estimates have been revisited and improved   by Tarr\`es in his PhD thesis (Tarr\`es 2000, 2001). 

The present section is clearly inspired by the work of these authors.  

\subsection{Mercer kernels}
\label{sec:mercer}
Recall that a Mercer kernel is a continuous symmetric function $K : M \times M \to \RR$ inducing a positive operator on $L^2(\lambda)$ (i.e., $\la Kf, f \ra_{\lambda} \geq 0).$ 
The following theorem is a fairly standard result in the theory of reproducing kernel Hilbert spaces (see e.g~ Aronszajn (1950) or Cucker and Smale (2001, Chapter III, 3)).
\bthm Let $K$ be a  Mercer kernel. Then there exists a unique Hilbert
space $\cH^K \subset C^0(M)$, the self reproducing space, such that
\bdes \iti For all $\mu\in\cM(M)$, $K\mu\in\cH^K$;
\itii For all $\mu$ and $\nu$ in $\cM(M)$,
\beq
\label{eqKnorm} \langle K\mu,K\nu\rangle_K = \int\int K(x,y)\mu(dx)\mu(dy). \eeq
\itiii $K(L^2(\lambda))$, $\{K_x,~x\in M\}$ and $K(\cM(M))$ are dense
in $\cH^K$.
\itiv For all $h\in\cH^K$ and $\mu\in\cM(M)$,
\beq \mu h = \langle K\mu,h\rangle_K.\eeq
\edes
Moreover, the mappings $K:\cM_s(M)\to \cH^K$ and $K:C^0(M)\to \cH^K$
are linear continuous and
for all $h\in \cH^K$,
\beq\label{ik} \|h\|_\infty\leq \|K\|_\infty^{1/2} \|h\|_K.\eeq
Hence,  the mapping $i_K:\cH^K\to C^0(M)$ defined by $i_K(h)=h$ is
continuous. \ethm


From now on and throughout the remainder of the section we  assume that hypothesis \ref{hyp:occ2} holds and we set 
\beq K = V_+ + V_- \eeq
where $V_+$ and $V_-$ have been defined by (\ref{defV+}).
According to hypothesis \ref{hyp:occ2}, $V_+$ and $V_-,$ hence  $K$ are Mercer kernels.
\bprop
\bdes
\iti One has the orthogonal decomposition (in $\cH^K$) $$\cH^K = \cH^{V_+}\oplus \cH^{V_-}.$$
\itii  Let  $\pi^+$ and  $\pi^-$ be the orthogonal projections onto
$\cH^{V_+}$ and onto $\cH^{V_-}$ (note that $\pi^{\pm}=\pi_{\pm}$
restricted to $\cH^{K}$). Then for all $h\in\cH^K$,
\beq \|h\|_K^2 = \|\pi^+h\|^2_{V_+} + \|\pi^-h\|^2_{V_-}. \eeq
\itiii  $V(\cM(M))=K(\cM(M))$ and  for all
$\mu\in\cM(M)$ and $h\in\cH^K$,
\beq\langle V\mu,h\rangle_K = \mu\pi^+h-\mu\pi^-h. \eeq
\edes
\eprop
\prf 
We have the
orthogonal decomposition (in $\cH^K$)
$K(L^2(\lambda))=V_+(L^2(\lambda))\oplus V_-(L^2(\lambda))$
(since $\langle V_+f,V_-g\rangle_K=\langle
K\pi_+f,K\pi_-g\rangle_K=\langle K\pi_+f,\pi_-g\rangle_\lambda=
0$). This implies  the orthogonal decomposition
$\cH^K = \cH^{V_+}\oplus \cH^{V_-},$
 because $\cH^{V_+}$ and $\cH^{V_-}$ are respectively the closures of
$V_+(L^2(\lambda))$ and of $V_-(L^2(\lambda))$ in $\cH^K$ (since
$\langle V_+f,V_+g\rangle_{V_+} = \langle V_+f,g\rangle_{\lambda} =
\langle K\pi_+f,\pi_+g\rangle_{\lambda} = \langle V_+f,V_+g\rangle_K)$.
Assertions $(ii)$ and $(iii)$ easily follow.
\qed


\brem \label{base}
{\rm
Let $(e_i)_i$ be an orthonormal basis of $\cH^K$ such that for all $i$,
$e_i$ belongs to $\cH^{V_+}$ or to $\cH^{V_-}$ and we set $\eps_i=\pm 1$ when
$e_i\in H^{V_\pm}$. Then we have
\beqarr
V_{\pm}(x,y) &=& \sum_i 1_{\eps_i=\pm 1} e_i(x) e_i(y), \\
K(x,y) &=& \sum_i e_i(x) e_i(y), \\
V(x,y) &=& \sum_i \eps_ie_i(x) e_i(y),
\eeqarr
the convergence being uniform by Mercer theorem (see e.g~ Chap XI-6 in Dieudonn\'e (1972) or
Cucker and Smale (2001)).}
\erem

\blem The mappings $V:\cM_s(M)\to \cH^K$ and $V:C^0(M)\to \cH^K$ are bounded operators. \elem
\prf This follows from the fact that for every $\mu\in\cM(M)$ and
every $f\in C^0(M)$
\beqarr
\|V\mu\|_{K}^2 &=& \mu^{\otimes 2}K \quad\leq\quad \|K\|_\infty \times |\mu|^2;\\ 
\|Vf\|_{K}^2 &\leq& \|K\|_\infty \times \|f\|^2_\infty. \qquad
\qed \eeqarr 

\subsection{The vector field $\tilde{Y} = \tilde{Y}_V$}
\label{sec:tidy}
We denote by $\cH^K_0$ the closure in $\cH^K$ of
$V(\cM_0(M))=K(\cM_0(M))$ and we set $\cH^K_1=V1+\cH^K_0$, the closure
of $V(\cM_1(M))=K(\cM_1(M))$. Equipped with the scalar product
$\langle\cdot,\cdot\rangle_{K}$, $\cH^K_0$ and $\cH^K_1$ are
respectively an Hilbert space and an affine Hilbert space.
 
We let $\tilde{Y} = \tilde{Y}_V:\cH^K_1\to \cH^K_0$ be the vector field defined by 
\beq
\label{deftilY}
 \tilde{Y}(h)=-h+V\xi(h).\eeq
Observe that $\tilde{Y}$ is exactly defined like the vector field $Y$ (introduced in the subsection \ref{sec:defY}) but for the fact that $\tilde{Y}$ is a vector field on  $\cH^K_1$ (rather than on $C^0(M)$).

Recall that we let $\Phi$ denote the smooth flow on $\cM_s(M)$ induced by the vector field $F$
defined in section \ref{sec:background} (equation (\ref{eq:ode})).
The proof of the  following lemma is similar to the proof of Lemma \ref{th:Ysink}.
\blem
\label{th:tilYsink}
The vector field $\tilde{Y}$ induces a global smooth flow $\tilde{\Phi}$
on $\cH^K_1(M).$ Furthermore
\bdes
\iti $V \Phi_t (\mu) = \tilde{\Phi}_t(V\mu)$ for all $\mu  \in \cM_s(M)$ and $t \in \RR.$
\itii $V$ maps homeomorphically $\mathsf{Fix}(\Pi)$ to $\tilde{Y}^{-1}(0),$  sinks to sinks and saddles to saddles.
\edes
\elem

\subsection{The stable manifold theorem and the function $\eta$}
\label{sec:stabman}
Let $\mu^*$ be a nondegenerate unstable fixed point of $\Pi$ and let 
\beq h^* = V\mu^*.\eeq
By Lemma \ref{th:tilYsink}, $h^*$ is a saddle for $\tilde{Y}$. Therefore
there exists constants $C$, $\lambda > 0$ and  a  splitting 
\beq\cH_0^K= H^s \oplus H^u,\eeq
with $H^u \neq \{0\},$
invariant under $D\tilde{\Phi}$
such that for all $t \geq 0$ and $u \in H^u$, 
\beqar
\label{eq:instable1}
\|D\tilde{\Phi}_t(h^*) u \|_{K} &\geq& C  e^{\lambda t} \|u\|_{K} \\
\mbox{and } \quad
\label{eq:instable2}
\|D\tilde{\Phi}_{-t}(h^*) u\|_{K} &\geq& C  e^{\lambda t} \|u\|_{K}.
\eeqar
\brem
{\rm 
Let, for $\alpha \in \RR,$ $H^\alpha = \{u\in L^2(\lambda), V\mathsf{T}(h^*)u=\alpha u\}$
where $\mathsf{T}(f)$ is the operator defined in  proposition \ref{th:spectrum}. From the proof of Lemma \ref{th:Ysink} it is easy to see
that
\beqarr
H^u &=& \oplus_{\alpha < - 1} H^{\alpha} \\
\mbox{and }\quad
H^s &=& \oplus_{\alpha > - 1} H^{\alpha}. \eeqarr
In particular, note that $H^u$ has finite dimension.
} 
\erem
\subsubsection*{The stable manifold theorem}
Set $h^* = h^*_s + h^*_u \in H^s \oplus H^u.$
By the {\em stable manifold theorem} (see e.g~ Hirsch and Pugh (1970) or Irwin (1970)) there exists
a neighborhood $\cN_0 = \cN_0^s \oplus \cN_0^u$ of $h^*$, with
 $\cN_0^s$ (respectively, $\cN_0^u$) a ball around $h^*_s$ in $H^s,$
(respectively, $h^*_u$ in $H^u$)
and a smooth function $\Gamma : \cN_0^s \to \cN_0^u$
such that
\bdes
\ita
$D\Gamma(h^*_s) = 0.$
\itb
The graph of $\Gamma:$
$$ Graph(\Gamma) = \{v + \Gamma(v)\, :  v \in \cN_0^s\},$$
equals the local stable manifold   of $h^*:$   
\beqarr
W^s_{loc}(h^*)  &=& \{h\in \cH^K_1: \quad 
   \forall t\geq 0,~\tilde{\Phi}_t(h) \in \cN_0 \\
&&\qquad \hbox{ and } \lim_{t\to\infty}\tilde{\Phi}_t(h)=h^*\}. \\
&=& \{h\in \cH^K_1:\quad \forall t\geq 0,~\tilde{\Phi}_t(h) \in
   \cN_0\}.
\eeqarr
\itc $W^s_{loc}(h^*)$ is an invariant manifold. That is for all $t \in \RR$,
$$\tilde{\Phi}_t(W^s_{loc}(h^*)) \cap \cN_0 \subset W^s_{loc}(h^*).$$
\edes
\subsubsection*{The function $\eta$}
Let $r : \cN_0 = \cN_0^s \oplus \cN_0^u \to W^s_{loc}(h^*)$
and $R :\cN_0 \to \RR$
be the functions defined by
$$r(h_s + h_u) = h_s + \Gamma(h_s)$$
and
$$R(h)=\|h-r(h)\|^2_{K}.$$ Then  $r$ and  $R$ are  smooth 
and $R$ vanishes  on  $W_{loc}^s(h^*)$.
\blem 
\label{th:lemR}
There exists $T>0$ and a neighborhood $\cN_1\subset \cN_0$ of $h^*$
in $\cH^K_1$ such that for all $h\in \cN_1$,
$\tilde{\Phi}_T(h)\in\cN_0$ and  
\beq R(\tilde{\Phi}_T(h))\geq R(h).\eeq
\elem  
\prf Using inequality (\ref{eq:instable1})   
we choose $T$ large enough so that for all 
$v\in H^u$, 
\beq\|D\tilde{\Phi}_T(h^*)v\|_{K}^2\geq 4 \|v\|^2_{K}.\eeq
Hence, there exists a neighborhood
$\cN_0'\subset \cN_0$ of $h^*$ such that for all $h\in\cN_0'$,
$\tilde{\Phi}_T(h)\in\cN_0,$ and for  all $v\in H^u$
\beq \|D\tilde{\Phi}_T(h)v\|^2_{K}\geq 3 \|v\|^2_{K}.\eeq
One may furthermore assume that
for all  $h\in\cN_0'$ (taking $\cN'_0$ small enough),
\beq ||D(r \circ \tilde{\Phi}_T)(h) - D(r \circ
\tilde{\Phi}_T)(h^*)||_{K} \leq 1.\eeq

Now, one has
\beq \label{eq:o}
\tilde{\Phi}_T(h) - \tilde{\Phi}_T(r(h)) -  D\tilde{\Phi}_T(r(h))(h-r(h)) = o(\|h-r(h)\|_{K}).\eeq
Using first the invariance of $W^s_{loc}(h^*)$, then equation
(\ref{eq:o}) with the fact that $D(r \circ \tilde{\Phi}_T)(h^*) v =
Dr(h^*) D \tilde{\Phi}_T(h^*) v = 0$ for all $v \in H^u$, we get
\beqarr
r(\tilde{\Phi}_T(h)) -  \tilde{\Phi}_T(r(h)) \quad =\quad  r(\tilde{\Phi}_T(h))
-  r (\tilde{\Phi}_T(r(h))) \hskip-180pt &&  \\
&=& D(r \circ \tilde{\Phi}_T)(r(h))(h-r(h)) + o(\|h-r(h)\|_{K})\\
&=& [ D(r \circ \tilde{\Phi}_T)(r(h)) - D(r \circ
  \tilde{\Phi}_T)(h^*)](h-r(h)) \\
&& \hskip100pt+\quad   o(\|h-r(h)\|_{K}). \eeqarr
Thus we obtain the upper-estimate
\beqarr
\|\tilde{\Phi}_T(h) - r(\tilde{\Phi}_T(h)) -
D\tilde{\Phi}_T(r(h))(h-r(h))\|_K \hskip-100pt &&\\
&\leq&   \|h-r(h)\|_{K} + o( \|h-r(h)\|_{K}).\eeqarr

This yields 
$$R(\tilde{\Phi}_T(h))\geq 2R(h) + o(R(h)).$$ 
We finish the proof of this lemma by taking $\cN_1\subset \cN_0$, a
neighborhood of $h^*$, such that for every $h\in \cN_1$, $o(R(h))\geq
-R(h)$. \qed
  
\medskip  
Let $\cN_2\subset \cN_1$ be a neighborhood of $h^*$  
such that for every $h\in \cN_2$ and every $t\in [0,T]$,  
$\tilde{\Phi}_{-t}(h)\in \cN_1$ ($T$ being the constant given in the previous
lemma). For every $h\in\cN_2$, set 
\beq \eta(h)=\int_0^T R(\tilde{\Phi}_{-s}(h))ds. \eeq
Then $\eta$ satisfies the following
\blem \label{lemeta} 
\bdes \iti $\eta(h)=0$ for every $h\in\cN_2\cap W^s_{loc}(h^*)$. 
\itii $\eta$ is $C^2$ on $\cN_2$.  
\itiii For every $h\in \cN_2$, $$D\eta(h)\tilde{Y}(h)\geq 0.$$ 
\itiv For every positive $\eps$ there exists $\cN^\eps_2\subset\cN_2$
and $D>0$ such that for all $h\in \cN_2^\eps$, $u$ and $v$ in
$\cH^K_0$,
\beqarr
|D^2_{u,v}\eta(h)-
D^2_{u,v}\eta(h^*)| &\leq& \eps \times \|u\|_{K}
\times \|v\|_{K}.\\
|D^2_{u,v}\eta(h^*)| &\leq& D \times \|u\|_{K}
\times \|v\|_{K}.
\eeqarr
\itv $D^2_{u,u}\eta(h^*)=0$ implies that $u\in
H^s.$
\itvi There exists a constant $C_\eta$ such that for all $u\in \cH^K_0$ and 
$h\in \cN_2$,  
$$|D\eta(h) u|\leq C_\eta\times\|u\|_{K} \times\sqrt{\eta(h)}.$$ 
\edes \elem  
\prf {\bf (i)} and {\bf (ii)} are clear. We have for $h\in\cN_2$
\beqarr  
D\eta(h)\tilde{Y}(h)  
&=& \lim_{s\to 0} \frac{1}{s}(\eta({\tilde{\Phi}}_s(h))-\eta(h)) \\  
&=& \lim_{s\to 0} \frac{1}{s}  
 \left(\int_0^s R({\tilde{\Phi}}_t(h))dt-\int_{T-s}^T R({\tilde{\Phi}}_{-t}(h))dt\right) \\  
&=& R(h)-R({\tilde{\Phi}}_{-T}(h)) \quad\geq\quad 0  \mbox{ (by Lemma \ref{th:lemR})}.
\eeqarr  
This shows {\bf (iii)}. Assertion {\bf (iv)} follows from the facts 
that $\eta$ is $C^2$.
For $h\in \cN_0$ and $u\in \cH^K_0$,
\beqarr 
DR(h)u &=& 2\langle h-r(h),u-Dr(h)u\rangle_{K}\\ 
D^2_{u,u}R(h) &=& 2\|u - Dr(h)u\|^2_{K} -  
2\langle h-r(h),D^2_{uu}r(h)\rangle_{K}.
\eeqarr 
Therefore
\beq D^2_{u,u}\eta(h^*) 
=2\int_0^T\|(I-Dr(h^*)) D\tilde{\Phi}_{-s}(h^*)u\|^2_{K} ds.\eeq
Since $Dr(h^*)$ is the projection onto $H^s$ parallel to $H^u$
one sees that $D^2_{u,u}\eta(h^*) = 0$ if and only if $D\tilde{\Phi}_{-s}(h^*) u \in H^u$ for all $s.$ 
This proves {\bf (v)} after remarking that for $s=0$,
$D\tilde{\Phi}_{-s}(h^*)u=u$.

We now prove {\bf (vi)}. For $u\in \cH^K_0$ and $h\in \cN_2$, 
$$D\eta(h)u = 2\int_0^T \langle h_s-r(h_s),u_s-Dr(h_s)u_s\rangle_{K} ds,$$ 
where $u_s=D\tilde{\Phi}_{-s}(h)u$ and $h_s=\tilde{\Phi}_{-s}(h)$. We
conclude using Cauchy-Schwartz inequality. \qed 

\subsection{Semigroups estimates} 
In the following, $\cD_2$ denotes the $L^2$-domain of 
the Laplacian on $M$. For $h\in C^1(M)$, 
set $A_h:\cD_2\to L^2(\lambda)$ defined by 
\beq A_h f=-\Delta f + \langle\nabla h,\nabla f\rangle,\eeq
and $Q_h:L^2(\lambda)\to\cD_2$ such that  
\beq 
\label{eq:defQ}
-Q_hA_hf=f-\langle\xi(h),f\rangle_\lambda.\eeq
Let $\P^h_t$ be the Markovian semigroup symmetric with respect to
$\mu_h=\xi(h)\lambda$ and with generator $A_h$. Note that $Q_h$ can be
defined by
\beq Q_hf = \int_0^\infty (\P^h_tf - \mu_h f)dt. \eeq
  
\blem \label{lem51} There exists a constant $K_1$ such that  
for all $f\in C^0(M)$ and $h\in \cH^K_1$ satisfying 
 $\|h\|_\infty\leq \|V\|_\infty$,  
$Q_h f\in C^1(M)\cap\cD_2$ and   
\beq\|\nabla Q_hf\|_\infty\leq K_1\|f\|_\infty.\eeq
\elem
\prf The proof of Lemma 5.1 in (BLR) can be easily adapted to prove this
lemma. \qed

\medskip
We denote by $C^{1,1}(M^2)$ the class of functions $f\in C^0(M^2)$ such
that for all $1\leq k,l\leq n$,
 $\frac{\partial}{\partial x^k}\frac{\partial}{\partial y^l}
f(x,y)$ exists and belongs to $C^0(M^2)$, where $(x^k)_k$ is a system of
local coordinates. For $f\in C^{1,1}(M^2)$, we define 
$\nabla^{\otimes 2} f\in C^0(TM\times TM)$ by
\beqarr \nabla^{\otimes 2} f((x,u),(y,v)) &=&
(\nabla_u\otimes \nabla_v) f(x,y) \\
&=& \sum_{k,l} u^kv^l \frac{\partial}{\partial x^k}
 \frac{\partial}{\partial y^l}f(x,y),
\eeqarr
in a system of local coordinates. We also define 
$\hbox{Tr}(\nabla^{\otimes 2}f)\in C^0(M)$, the trace of
$\nabla^{\otimes 2}f$, by
($d$ denotes the dimension of $M$)
$$\hbox{Tr}(\nabla^{\otimes 2}f)(x) = \sum_{k=1}^d
\frac{\partial}{\partial x^k}  \frac{\partial}{\partial y^k}f(x,x).$$
This definition is of course independent of the chosen system of local
coordinates.

\brem\label{rem51} Lemma \ref{lem51} implies that for all $f\in C^0(M^2)$ and $h\in \cH^K_1$ satisfying 
 $\|h\|_\infty\leq \|V\|_\infty$,
$Q_h^{\otimes 2} f\in C^{1,1}(M^2)$ and
\beq\|\nabla^{\otimes 2} Q_h^{\otimes 2}f\|_\infty\leq
K_1^2\|f\|_\infty.\eeq
This also implies that
\beq\|\hbox{Tr}(\nabla^{\otimes 2}Q_h^{\otimes 2}f)\|_\infty\leq dK_1^2\|f\|_\infty.\eeq
\erem
 
\blem \label{lemlip} There exists a constant $K_2(=K_1^2)$ such that for 
all $f\in C^0(M)$, $h_1$ and $h_2$ in $\cH^K_1$ satisfying 
$\|h_1\|_\infty\vee\|h_2\|_\infty\leq \|V\|_\infty$, 
\beq \|\nabla Q_{h_2}f-\nabla Q_{h_1}f\|_\infty 
\leq K_2\|f\|_\infty\|\nabla h_2-\nabla h_1\|_\infty.\eeq
\elem  
\prf Set $u=Q_{h_1}f$. Then  
$$-A_{h_1}u=f-\langle\xi(h_1),f\rangle_\lambda$$  
and since $A_{h_2}u-A_{h_1}u=\langle\nabla (h_2-h_1),\nabla u\rangle$,
\beqarr  
Q_{h_2}f  
&=& -Q_{h_2}(A_{h_1}u-\langle\xi(h_1),f\rangle_\lambda)\\  
&=& -Q_{h_2}A_{h_1}u \\  
&=& -Q_{h_2}A_{h_2}u + Q_{h_2}f_h \eeqarr  
where $h=h_2-h_1$ and $f_h=\langle\nabla h,\nabla u\rangle$. Thus  
$$Q_{h_2}f=Q_{h_1}f-\langle\xi(h_2),Q_{h_1}f\rangle_\lambda  
+ Q_{h_2}f_h$$  
and  
$$\nabla Q_{h_2}f-\nabla Q_{h_1}f = \nabla Q_{h_2}f_h.$$ 
Lemma \ref{lem51} implies that 
\beqarr \|\nabla Q_{h_2}f_h\|_\infty &\leq& K_1\|f_h\|_\infty\\ 
\mbox{and }\quad \|\nabla Q_{h_1}f\|_\infty &\leq& K_1\|f\|_\infty. \eeqarr 
We conclude since  
$\|f_h\|_\infty\leq \|\nabla h\|_\infty \|\nabla Q_{h_1}f\|_\infty$. 
\qed 

\brem \label{remlip} Lemma \ref{lemlip} implies that for all $f\in
C^0(M^2)$, $h_1$ and $h_2$ in $\cH^K_1$ satisfying
$\|h_1\|_\infty\vee\|h_2\|_\infty\leq \|V\|_\infty$, 
\beq\|\nabla^{\otimes 2} Q_{h_2}^{\otimes 2}f - \nabla^{\otimes 2}
Q_{h_1}^{\otimes 2}f\|_\infty  \leq K_2^2\|f\|_\infty\|\nabla
h_2-\nabla h_1\|_\infty^2.\eeq
This implies that
\beq\|\hbox{Tr}(\nabla^{\otimes 2} (Q_{h_2}^{\otimes 2} -
Q_{h_1}^{\otimes 2})f)\|_\infty  \leq dK_2^2\|f\|_\infty\|\nabla
h_2-\nabla h_1\|_\infty^2.\eeq
\erem

\subsection{It\^o calculus}
\label{sec:itocal} 
Set $h_t=V\mu_t.$
Given a smooth (at least $C^2$) function 
$$\begin{array}{lll}
\RR \times M &\to& \RR \\
(t,x) &\mapsto&  F_t(x),\end{array}$$
It\^o's formula reads
\beq
\label{eq:ito}
dF_t(X_t) = \partial_t F_t(X_t) dt + A_{h_t} F_t(X_t) dt + dM_t
\eeq
where $M$ is a martingale with ($\langle\cdot,\cdot\rangle_t$ denotes
the martingale bracket)
$$ \frac{d}{dt}\langle M^f\rangle_t =   
\frac{1}{t^2}\|\nabla F_t(X_t)\|^2. $$
  
Set  $Q_t=Q_{h_t}$ and $F_t(x) =  \frac{1}{t}Q_tf(x)$ for some $f \in C^0(M).$
Then (\ref{eq:ito}) (note that It\^o's formula also holds if
$(t,x)\mapsto F_t(x)$ is $C^1$ in $t$ and  for all
$t$, $F_t\in\cD_2$, which holds here) combined with (\ref{eq:defQ}) 
 gives
\beq
\label{eq:ito2}
d\left(\frac{1}{t}Q_tf(X_t)\right)  = 
 \frac{H_t f}{t^2}dt +  \frac{\la \xi(h_t),f\ra_{\lambda} - f(X_t)}{t}  + dM^f_t \eeq  
where $H_t$ is the measure defined by  
\beq H_t f = - Q_tf(X_t) + t\left(\frac{d}{dt}Q_t\right)f(X_t), \eeq  
$M^f$ is a martingale with  
\beq \frac{d}{dt}\langle M^f\rangle_t =   
\frac{1}{t^2}\|\nabla Q_t f(X_t)\|^2. \eeq  
Using the fact that
$$\frac{d}{dt}\mu_t f = \frac{f(X_t) - \mu_t f}{t}$$
together with the definition of the vector field $F,$ (\ref{eq:ito2})
can be rewritten as (recall that $F(\mu)=-\mu+\Pi(\mu)$ and that
$\Pi(\mu)=\xi(V\mu)\lambda$)
\beq 
\label{eq:ito3}
d\mu_t f = \frac{F(\mu_t) f}{t}dt  
-d\left(\frac{1}{t}Q_tf(X_t)\right) +  \frac{H_t f}{t^2}dt + dM^f_t \eeq  
Note that there exists a constant $H$ such that for all $t\geq 0$ and 
$f\in C^0(M)$,   $|H_tf|\leq H\|f\|_\infty$ (see Lemmas 5.1 and 5.6 in (BLR)).  
\smallskip

Let $\nu_t$ be the measure defined by  
\beq \nu_t f = \mu_t f +\frac{1}{t}Q_tf(X_t),\qquad f\in C^0(M).\eeq  
Then $|\mu_t - \nu_t|\to 0$ and   
\beq d\nu_t f = \frac{F(\nu_t) f}{t}dt +  
\frac{N_t f}{t^2}dt + dM^f_t, \eeq  
with $N_t$ the measure defined by
$N_tf=H_tf+t\left(F(\mu_t)-F(\nu_t)\right)f$. Since $F$ is Lipschitz,
there exists a constant $N$ such that for all $t\geq 0$ and $f\in
C^0(M)$, 
\beq
\label{eq:boundN}
|N_tf|\leq N\|f\|_\infty.
\eeq

For every $t\geq 1$, set $g_t=V\nu_t$. Then using the fact that 
$VF(\mu)=\tilde{Y}(V\mu)$, 
\beq \label{dg} dg_t(x) = \frac{\tilde{Y}(g_t)(x)}{t}dt +  
\frac{N_tV_x}{t^2}dt + dM^{V_x}_t, \eeq  
where $V_x(y)=V(x,y)$. 

Note that $(g_t)_{t\geq 1}$ is a $\cH^K_0$-valued continuous
semimartingale. We denote its martingale part $M_t$, with
$M_t(x)=M^{V_x}_t-M^{V_x}_1$. In the following, $(e_i)$ denotes an
orthonormal basis of $\cH^K$ like in remark \ref{base}. Then
$M_t=\sum_i M^i_t e_i$, with $M^i_t=\langle M_t,e_i\rangle_{K}$.
Using the fact that for all $\mu\in \cM_0(M)$, 
$$\langle M_t,K\mu\rangle_K=\int M_t^{V_x}\mu(dx)$$
we have
\beqarr
\frac{d}{ds}\langle \langle M_\cdot ,K\mu\rangle_{K}\rangle_s
&=& \int\int \frac{d}{ds}\langle M^{V_x},M^{V_y}\rangle_s
\mu(dx)\mu(dy) \\
&=&\int\int \frac{1}{s^2}\times \langle \nabla Q_sV_x(X_s), 
\nabla Q_sV_y(X_s) \rangle \mu(dx)\mu(dy) \\
&=&\frac{1}{s^2}\times \| \nabla Q_sV\mu(X_s) \|^2.
\eeqarr

This implies that for $h$ in $\cH^{V_+}$ or in $\cH^{V_-}$
\beq  \label{varquad1}
\frac{d}{ds}\langle \langle M_\cdot ,h\rangle_{K}\rangle_s
= \frac{1}{s^2}\times \| \nabla Q_sh(X_s) \|^2 \eeq
and
\beq  \label{varquad2}
\frac{d}{ds}\langle M^i,M^j \rangle_s
= \frac{\eps_i\eps_j}{s^2}\times \langle \nabla Q_s e_i(X_s), \nabla Q_s e_j(X_s)
\rangle. \eeq
\blem \label{martl2}
There exists a constant $C_1$ such that for every $s\geq 1$,
\beq\E[\|M_s\|_{K}^2]\leq C_1.\eeq
\elem
\prf We have
\beqarr
\frac{d}{ds} \E[\|M_s\|_{K}^2]
&=& \sum_i \frac{d}{ds} \E[\langle M^i,M^i \rangle_s]\\
&=& \frac{1}{s^2}\times \E\left[\sum_i \|\nabla Q_s e_i(X_s)\|^2\right]\\
&=& \frac{1}{s^2}\times \E\left[\hbox{Tr}(\nabla^{\otimes 2}
  Q_s^{\otimes 2} K) (X_s,X_s)\right]
\eeqarr
since $K=\sum_i e_i\otimes e_i$. We conclude using remark
\ref{rem51} and taking $C_1=dK_1^2\|K\|_\infty$. \qed
\subsection{A first lemma}
\label{sec:lemmafirst}
Let $L$ be a positive constant we will fix later on. 
Set $\eta_t=\eta(g_t)1_{g_t\in\cN_2}$ where $\cN_2$ is like in Lemma \ref{lemeta}. Let $\cN$ be a neighborhood of $\mu^*$
(for the narrow topology).
For every $t\geq 1$, set $S_t=\inf\{s>t,~ \eta_s\geq L^2/s\}$ and  
$U^\cN_t=\inf\{s>t,~\mu_s\not\in \cN\}$. The purpose of this section is to
prove
\blem \label{lem1}  
There exist a neighborhood $\cN$ of $\mu^*$, $p\in ]0,1]$ and  
$T_1>0$ such that for all $t>T_1$,   
\beq \P[S_t\wedge U^\cN_t<\infty|\cB_t] \geq p.\eeq 
where $\cB_t$ is the sigma field generated by $\{B_s^i: \:  i = 1 \ldots N, s \leq t \}.$ \elem  
\prf We fix $\eps>0$. 
Since $V:\cP_w(M)\to \cH^K$ is continuous and $|\nu_t - \mu_t| \rar 0$ there exist $\tau_1$ large
enough and $\cN_\eps$ a neighborhood of $\mu^*$ such that for all
$t\geq \tau_1$, $\mu_t\in \cN_\eps$ implies that $\nu_t\in
V^{-1}(\cN_2^\eps)$, where $\cN_2^\eps$ is the neighborhood defined in
lemma \ref{lemeta}.
In particular, $\mu_t\in\cN_\eps$ implies that
$g_t=V\nu_t\in\cN_2^\eps$.

For every neighborhood $\cN\subset \cN_\eps$ of $\mu^*$ and every
$s\in [t,U^\cN_t]$, $ \eta_s=\eta(g_s)$. Then It\^o's formula with
formulas (\ref{dg}) and (\ref{varquad2}) gives
\beqar  
d\eta(g_s) &=&
\frac{D\eta(g_s)\tilde{Y}(g_s)}{s}ds +  
\frac{D\eta(g_s)(VN_s)}{s^2}ds + dM^\eta_s\label{ItoH} \\  
&& +\quad \frac{1}{2} \sum_{i,j} D^2_{i,j}\eta(g_s)
\times\langle \eps_i\nabla Q_se_i(X_s) , 
\eps_j\nabla Q_se_j(X_s)\rangle\times \frac{ds}{s^2},\non
\eeqar
where $VN_s(x)=N_sV_x$ and $M^\eta$ is the martingale defined by
\beq \label{defmeta} dM^\eta_s=D\eta(g_s)dM_s.\eeq

We now intend to prove that
\beq \E[\eta(g_{S_t\wedge U^\cN_t})|\cB_t]-\eta(g_t) \geq  
-C\eps/t+(K^*/t)\P[S_t\wedge U^\cN_t=\infty|\cB_t],\eeq
where $C$ and $K^*$ are positive constants. In order to do this, we
bound from below the four terms in the right hand side of (\ref{ItoH}).

Lemma \ref{lemeta} {\bf (iii)} implies that
$D\eta(g_s)\tilde{Y}(g_s)\geq 0$.
Using Lemma \ref{lemeta} {\bf (vi)} and inequality (\ref{eq:boundN}), it can be easily seen that there 
exists a constant $N_\eta$ such that for $s\in [t,U_t^\cN]$ 
$$|D\eta(g_s)VN_s|\leq N_\eta\sqrt{\eta(g_s)}.$$
Then 
$$\int_t^{S_t\wedge U^\cN_t}\frac{D\eta(g_s)VN_s}{s^2}ds  
\geq -LN_\eta\int_t^\infty \frac{ds}{s^{5/2}}.$$  
We choose $\tau_2\geq \tau_1$ large enough such that for all $t\geq
\tau_2$,
\beq \label{eqln} LN_\eta\int_t^\infty \frac{ds}{s^{5/2}}\leq \eps/t. \eeq
This gives an estimate of the second term. Since the third term is a
martingale increment, after taking the expectation, this term will
vanish.

We now estimate the last term.
For $s>0$, set \beq\label{defAs}
\Gamma_s =  \sum_{i,j} D^2_{i,j}\eta(g_s)
\times\langle \eps_i\nabla Q_se_i(X_s) , 
\eps_j\nabla Q_se_j(X_s)\rangle \eeq
and, for $\mu\in\cP(M)$ and $x\in M$, set
\beq\label{defA}\Gamma(\mu,x) =  \sum_{i,j} D^2_{i,j}\eta(h^*)
\times\langle \eps_i\nabla Q_{V\mu}e_i(x) , 
\eps_j\nabla Q_{V\mu}e_j(x)\rangle.\eeq
Lemma \ref{lemeta} {\bf (iv)} implies that for $s\in [t,U_t^\cN]$
(to prove this upper-estimate, one can use a system of local
coordinates and use the fact that $K=\sum_i e_i\otimes e_i$)
\beqarr
|\Gamma_s-\Gamma(\mu_s,X_s)| 
&\leq& \eps\times\sum_i\|\nabla Q_se_i(X_s)\|^2\\
&\leq& \eps\times \hbox{Tr}(\nabla^{\otimes 2}Q_s^{\otimes 2} K)(X_s)
\eeqarr
Thus $|\Gamma_s-\Gamma(\mu_s,X_s)|\leq C_1\times\eps$ where $C_1$ is the same
constant as the one given in Lemma \ref{martl2}.

\blem \label{Acont} $\Gamma:\cP_w(M)\times M\to \RR^+$ is continuous. \elem
\prf We only prove the continuity in $\mu$.
For $\mu$ and $\nu$ in $\cP(M)$ and
$x\in M$,
$$\Gamma(\mu,x)-\Gamma(\nu,x)= \sum_{i,j} D^2_{i,j}\eta(h^*)
\langle u_i(\mu,x) -  u_i(\nu,x),  u_i(\mu,x) +  u_i(\nu,x) \rangle$$
where $u_i(\mu,x)=\eps_i\nabla Q_{V\mu}e_i(x)$. Using lemma
\ref{lemeta} {\bf (iv)} and Cauchy-Schwartz inequality,
\beqarr|\Gamma(\mu,x)-\Gamma(\nu,x)| &\leq& D\times 
\left(\hbox{Tr}(\nabla^{\otimes 2} (Q_{V\mu}-Q_{V_\nu})K)(x)\right)^{1/2}\\
&&\times\left(\hbox{Tr}(\nabla^{\otimes 2} (Q_{V\mu}+Q_{V_\nu})K)(x)\right)^{1/2}.\eeqarr
Remarks \ref{rem51} and \ref{remlip} imply that
$$ |\Gamma(\mu,x)-\Gamma(\nu,x)| \leq D\times \sqrt{2}dK_2K_1\|K\|_\infty \times \|\nabla V\mu - \nabla
V\nu\|_\infty$$
which converges towards $0$ as $\hbox{dist}_w(\mu,\nu)\to 0$. The
proof of the continuity in $x$ is similar. \qed

\medskip
Lemma \ref{Acont} implies that we can choose the neighborhood
$\cN\subset\cN_\eps$ of $\mu^*$ such that for all $s\in [t,U_t^\cN]$,
\beq |\Gamma(\mu_s,X_s)-\Gamma(\mu^*,X_s)|\leq \eps. \eeq
We now set $\Gamma^*(x)=\Gamma(\mu^*,x)$. Thus we now have
\beqar \Gamma_s &=& (\Gamma_s - \Gamma(\mu_s,X_s)) +
(\Gamma(\mu_s,X_s))-\Gamma^*(X_s)) + \Gamma^*(X_s)\non\\
&\geq& -(C_1+1)\times\eps + \Gamma^*(X_s). \label{eq:gamma*} \eeqar

\medskip 
Finally using (\ref{eqln}) and (\ref{eq:gamma*}) (with the convention $ \eta_{S_t\wedge U^\cN_t}=0$ when  
$S_t\wedge U^\cN_t=\infty$)  
\beqarr  
\E[ \eta_{S_t\wedge U^\cN_t}|\cB_t]- \eta_t  
&\geq& -\frac{(2+C_1)\eps}{t}\\  
&& +\quad \frac{1}{2}\E\left[\left.\int_t^{\infty}\Gamma^*(X_s)\frac{ds}{s^{2}}  
                      1_{\{S_t\wedge U^\cN_t=\infty\}}\right|\cB_t\right] . 
\eeqarr

For all $s$, set $K(s)=\mu_s\Gamma^*$. 
Since $\Gamma^*(X_s)=K(s)+sK'(s)$ (recall that
$\mu_s=\frac{1}{s}\int_0^s\delta_{X_u}du$), integrating by
parts we get
$$\int_t^{\infty} \Gamma^*(X_s)\frac{ds}{s^{2}} = -\frac{K(t)}{t} +  
2\int_t^\infty \frac{K(s)}{s^2}ds.$$  
   
Since $\mu\mapsto \mu \Gamma^*$ is continuous, we can 
choose the neighborhood $\cN$ of $\mu^*$ such that for all $\mu\in 
\cN$, 
$$|\mu \Gamma^*-K^*|<\eps/3,$$
where $K^*=\mu^*\Gamma^*$.  
Then, on the event $\{S_t\wedge U^\cN_t=\infty\}$, for all  
$s\geq t$,  $$|K(s)-K^*|<\eps/3$$  
and
$$\int_t^{\infty} \Gamma^*(X_s)\frac{ds}{s^{2}} \geq \frac{K^*-\eps}{t}.$$  
Thus,  
\beq \E[ \eta_{S_t\wedge U^\cN_t}|\cB_t]- \eta_t \geq  
-(3+C_1)\eps/t+(K^*/t)\P[S_t\wedge U^\cN_t=\infty|\cB_t].\eeq  
\blem The constant $K^*=\int \Gamma^*(x)\mu^*(dx)$ is positive. \elem 
\prf We first remark that for all $f$ and $g$ in $C^0(M)$,
\beqarr
\langle \nabla Q_{h^*}f,\nabla Q_{h^*}g\rangle_{\mu^*}
&=& \langle f-\mu^*f,Q_{h^*}g\rangle_{\mu^*}\\
&=& \int_0^\infty \langle f-\mu^*f,
\P^{h^*}_t(g-\mu^*g)\rangle_{\mu^*} dt\\
&=& \int_0^\infty \langle \P^{h^*}_{t/2}(f-\mu^*f) , 
\P^{h^*}_{t/2}(g-\mu^*g)\rangle_{\mu^*} dt.
\eeqarr
Using this relation we get that
\beqarr 
K^* &=& 
 \sum_{i,j} D^2_{i,j}\eta(h^*)
\times \langle\eps_i\nabla Q_{h^*}e_i , \eps_j\nabla Q_{h^*}e_j\rangle_{\mu^*}\\
&=& \int_0^\infty \sum_{i,j} D^2_{i,j}\eta(h^*)
\times \langle \eps_i(\P_{t/2}^{h^*}e_i -\mu^*e_i)  , 
\eps_j(\P_{t/2}^{h^*}e_j-\mu^*e_j)\rangle_{\mu^*} dt\\
&=&  \int_0^\infty\int D^2\eta(h^*)(u_t^x,u_t^x) ~\mu^*(dx)\times dt,
\eeqarr 
where 
\beqarr
u_t^x &=& \sum_i \eps_i(\P_{t/2}^{h^*}e_i(x) -\mu^*e_i)e_i\\
&=& V(\P_{t/2}^{h^*}(x))-V\mu^* \eeqarr
($\P_{t/2}^{h^*}(x)$ denotes the measure
defined by $\P_{t/2}^{h^*}(x)f=\P_{t/2}^{h^*}f(x)$).

If $K^*=0$, then for all $x\in M$ and
$t\geq 0$, $u_t^x\in H^s$
since $D^2_{u,u}\eta(h^*)=0$ implies $u\in H^s$.
Thus, for all $x\in M$, $V_x-V\mu^*\in H^s,$ and for all $x$ and
$y$ in $M$, $V_x-V_y\in H^s$.
Therefore for every $\mu\in\cM_0(M)$, 
$V\mu\in H^s.$
This proves that $\cH^K_0\subset H^s$ and $H^u=\{0\}$. This gives 
a contradiction since the dimension of $H^u$ is larger than $1$. \qed
 
\medskip
On the other hand,  
$$\E[ \eta_{S_t\wedge U^\cN_t}|\cB_t]- \eta_t\leq  
\E[L^2/S_t\wedge U^\cN_t|\cB_t].$$ 
Therefore  
\beq L^2 \E[t/S_t\wedge U^\cN_t|\cB_t]\geq -(3+C_1)\eps +  
K^*\P[S_t\wedge U^\cN_t=\infty|\cB_t],\eeq  
and, since  
$$\P[S_t\wedge U^\cN_t<\infty|\cB_t]\geq \E[t/S_t\wedge 
U^\cN_t|\cB_t],$$ 
we have  
\beq \P[S_t\wedge U^\cN_t<\infty|\cB_t]\geq \frac{K^*-(3+C_1)\eps}{L^2+K^*}.\eeq  
Choosing $\eps<K^*/(3+C_1)$, this proves the lemma. \qed

\subsection{A second lemma}
\label{sec:lemmasec} 
We choose $\cN$, $p$ and $T_1$ like in lemma \ref{lem1}. Set  
\beq H=\{\liminf  \eta_t>0\}.\eeq
\blem \label{lem2}  
There exists $T_2>0$ such that for all $t>T_2$, on the event 
$\{S_t<U^\cN_t\}$,   
\beq\P[H|\cB_{S_t}]\geq 1/2.\eeq \elem  
  
\prf Fix $t>0$. Set  
\beq I_t=\inf_{s\in [S_t,U^\cN_t]} (M^\eta_s-M^\eta_{S_t}) \eeq
and  
\beq T_t=\inf\{s>S_t,  \eta_s=0\}. \eeq
  
On the event $\{S_t<U^\cN_t\}\bigcap\{I_t\geq -\frac{L}{2\sqrt{S_t}}\}$,  
for $s\in [S_t,T_t\wedge U^\cN_t]$, for some constant $N'<\infty$ we have  
\beqarr  
 \eta_s &=&  \eta_{S_t} + \int_{S_t}^s D\eta(g_u)\tilde{Y}(g_u)\frac{du}{u} +  
\int_{S_t}^s D\eta(g_u)VN_u \frac{du}{u^2} + M^\eta_s-M^\eta_{S_t} \\  
&\geq& \frac{L}{\sqrt{S_t}} - \frac{N'}{S_t} -  
\frac{L}{2\sqrt{S_t}}  
\quad\geq\quad \frac{L}{4\sqrt{S_t}}\eeqarr  
for $t\geq T_2$ large enough. Thus, for $t\geq T_2$,  
$$\liminf_{s\to\infty}  \eta_s\geq\frac{L}{4\sqrt{S_t}}$$  
and  
$$\{S_t<U^\cN_t\} \bigcap 
\left\{I_t\geq -\frac{L}{2\sqrt{S_t}}\right\}  
\subset H.$$  
  
Now, on the event $\{S_t<\infty\}$,  
\beqarr  
\P\left[\left. I_t < -\frac{L}{2\sqrt{S_t}}  
                 \right|\cB_{S_t}\right]  
&=& \P\left[\left.  
       \sup_{s\in [S_t,U^\cN_t]} -(M^\eta_s-M^\eta_{S_t})  
          > \frac{L}{2\sqrt{S_t}}  
                 \right|\cB_{S_t}\right]\\  
&\leq& \frac{4S_t}{L^2}\times\E\left[\left.\langle M^\eta\rangle_{U^\cN_t} 
           -\langle M^\eta\rangle_{S_t} \right| \cB_{S_t}\right]  
\eeqarr  
by Doob inequality. For $s\in [S_t,U^\cN_t]$,
\beqarr
d\langle M^\eta\rangle_s
&=& \sum_{i,j} D_i\eta(g_s) D_j\eta(g_s)d\langle M^i,M^j\rangle_s\\
&=& \frac{ds}{s^2}
\sum_{i,j} D_i\eta(g_s) D_j\eta(g_s)
\langle\eps_i \nabla Q_se_i(X_s),\eps_j\nabla Q_se_j(X_s)\rangle_s.
\eeqarr
Lemma \ref{lemeta} {\bf (vi)} implies that (recall that $K=\sum_i
e_i\otimes e_i$)
$$ \frac{d}{ds}\langle M^\eta\rangle_s
\leq \frac{1}{s^2} C_\eta^2\times
\hbox{Tr}(\nabla^{\otimes 2}Q_s^{\otimes 2} K)(X_s)
\leq  \frac{C}{s^2} $$
with  $C=C_1C_\eta^2$. Thus $\langle M^\eta\rangle_{U^\cN_t}-\langle
M^\eta\rangle_{S_t} \leq C/S_t$ and on the event  
$\{S_t<\infty\}$, we have  
$$\P\left[\left. I_t < -\frac{L}{2\sqrt{S_t}}  
                 \right|\cB_{S_t}\right]  
                     \leq 4C/L^2 .$$  
We choose $L$ such that $4C/L^2<1/2$. Then for $t\geq T_2$, on the
                 event $\{S_t<U^\cN_t\}$, 
$$  \P[H|\cB_{S_t}] \geq   
 \P\left[\left. I_t \geq -\frac{L}{2\sqrt{S_t}}  
                 \right|\cB_{S_t}\right] \geq 1/2. $$  
This proves the lemma. \qed  
  
\subsection{Proof of Theorem \ref{th:main3}}
We fix $\cN,p,T_1$ and $T_2$ like in lemmas \ref{lem1} and  
\ref{lem2}. Let $A=\{\exists t,~U^\cN_t=\infty\}$. Then for $t\geq  
T=T_1\vee T_2$, using lemmas \ref{lem1} and \ref{lem2},  
\beqarr  
\P[H|\cB_{t}]  
&\geq& \E[1_H1_{S_t<U^\cN_t}|\cB_{t}]\\  
&\geq& \E\left[\left.  
      \P[H|\cB_{S_t}]1_{S_t<U^\cN_t}  
                   \right|\cB_{t}\right]\\  
&\geq& \frac{1}{2}\times \P[S_t<U^\cN_t|\cB_t]\\  
&\geq& \frac{1}{2}\left(p-\P[U^\cN_t<\infty|\cB_t]\right). \eeqarr  
On one hand,  
$$\lim_{t\to\infty}\P[H|\cB_t]=1_H, \qquad\hbox{a.s.}$$  
On the other hand,  
$$\lim_{t\to\infty}1_{\{U^\cN_t=\infty\}}=1_A \qquad\hbox{a.s.}$$  
and  
\beqarr  
\E[|1_A-\P[U^\cN_t=\infty|\cB_t]|]  
&\leq& \E[|1_A-\P[A|\cB_t]|] \\  
&&+\quad  
\E[|\P[A|\cB_t]-\P[U^\cN_t=\infty|\cB_t]|]\\  
&\leq& \E[|1_A-\P[A|\cB_t]|] \\  
&&+\quad  
\E[|1_A-1_{\{U^\cN_t=\infty\}}|],  
\eeqarr  
which converges towards $0$ as $t\to\infty$. Thus  
$\lim_{t\to\infty}\P[U^\cN_t<\infty|\cB_t]=1_{A^c}$ in $L^1$ and  
\beq \label{eqH} 1_H\geq \frac{1}{2}(p-1_{A^c})\qquad \hbox{a.s.}\eeq  
This implies that a.s., $A\subset H$. But since  
$H\subset\{\mu_t\not\to\mu^*\}$ and $\{\mu_t\to\mu^*\}\subset A$, we  
have $\{\mu_t\to\mu^*\}\subset\{\mu_t\not\to\mu^*\}$ a.s. This implies  
that $\P[\mu_t\to\mu^*]=0$. \qed


\section{Appendix}
\label{sec:appendix}
Recall that we let 
 ${\cal G}$ denote the  set of $V \in C^k_{sym}(M \times M)$ such that $\Pi_V$ has nondegenerate fixed points.  Our purpose here is to prove Theorem
 \ref{th:generic}. That is that  ${\cal G}$ is open and dense.

\medskip
\paragraph{Openess.} 
We first prove that ${\cal G}$ is open. 
Let $V^* \in {\cal G}.$ Then the zeros of  $X_{V^*}$  are isolated (by
the inverse function theorem) and since $(X_{V^*})^{-1}(0)$ is compact
(Lemma \ref{th:fixzero})  ${X_{V^*}}^{-1}(0)$ is a finite set. Say
${X_{V^*}}^{-1}(0) = \{f_1, \ldots, f_d\}.$  

By the implicit function theorem applied to the map $(V,f) \mapsto X_V(f)$, there exist open  neighborhoods $U_i$ of $f_i,$ $W_i$ of $V^*$ 
and smooth maps $R_i : W_i \to U_i$ such that
\bdes
\ita $X_V(f) = 0 \Leftrightarrow f = R_i(V),$
for all $V \in W_i, f \in U_i,$ 
\itb $R_i(V^*) = f_i,$
\itc $DX_V(f)$ is invertible at $f = R_i(V).$
\edes
 It remains to show  that there exists an open neigborhood of $V^*$  $W \subset \bigcap_i W_i$ such that for all $V \in W$  equilibria of $X_V$ lie in $\bigcup U_i.$ In view of $(a)$ and $(c)$ above this will imply that $W \subset {\cal G}$ concluding the proof of openess.
Assume to the contrary that there is no such neighborhood. Then there exists $V_n \to V^*$ and $f_n \in {\cal B}_1 \setminus \bigcup_i U_i$ such that
$X_{V_n}(f_n) = 0.$
That is  \beq
\label{eq:fnfix}
f_n = \xi(V_n f_n) 
\eeq
Then by Lemma \ref{th:Vcompact},  we can extract from   $\{V^* f_n\}$ a subsequence $\{V^* f_{n_k}\}$  converging to some $g \in C^0(M).$ 
Now, $||V_n f_n - V f_n||_{\infty} \leq ||V_n - V^*||_{\infty}.$ Thus 
$V_{n_k} f_{n_k} \to g.$ Equation (\ref{eq:fnfix}) then implies  that $f_{n_k} \to f = \xi(g)$ and  $f = \xi(V^*f).$  Hence $f \in \bigcup_i U_i.$ A contradiction.
\paragraph{Density.}
We now pass to the proof of the density. 
Recall that if  $Z$ is a smooth map from one  Banach manifold to another,  a point  $h \in {\cal B}_2$  is called a {\em regular value} of $Z$ provided $DZ(f)$ is surjective for all $f \in Z^{-1}(h).$
Here, saying that $0$ is a regular value for $X_V$ is equivalent to saying that $X_V$ has nondegenerate equilibria.

Let ${\cal B}_1^k = {\cal B}_1 \cap C^k(M), {\cal B}_0^k = {\cal B}_0 \cap C^k(M)$ and
${{\cal B}^{+,k}_1} = {{\cal B}^+_1} \cap C^k(M).$ For all $V \in C^k_{sym}(M \times M)$ let   $Z_V : {{\cal B}^{+,k}_1} \to {\cal B}_0^k$ denote the $C^{\infty}$ vector field  defined by 
$$Z_V(f) = Vf + \log(f) - <Vf + \log(f), \one >.$$
Remark that for all $h \in {\cal B}_0^k$
$$DJ_V(f) h = <Z_V(f), h>.$$
Hence, by Proposition \ref{th:spectrum}, $X_V$ and $Z_V$ have the same set of equilibria
and $0$ is a regular value for $X_V$ if and only if it is a regular value for $Z_V.$ 

Given $h \in {\cal B}_0^k$ Let $V[h]$ be the symmetric function defined
by 
$$V[h](x,y) = V(x,y) - h(x) - h(y).$$
One has
$$Z_{V[h]}(f) = Z_V(f)  - h.$$
Therefore, $h$ is a regular value of $Z_V$ if and only if $0$ is a regular value of $Z_{V[h]}$ or, equivalently, a regular value of $X_{V[h]}.$

We claim that $Z_V$ is a {\em Fredholm map.} That is, a map whose
 derivative $D Z_V(f)$  is a Fredholm operator  for each  $f \in
 {{\cal B}^{+,k}_1}$ (see Section \ref{sec:proof1} for the definition
 of a Fredholm operator).   Hence by a theorem of Smale (1965)
 generalyzing Sard's theorem to Fredholm maps)  $\mathsf{R}_{Z_V}$
 is a residual (i.e., a countable intersection of open dense sets) set.
Being residual, it is dense. 
Therefore, for any $\eps > 0$ we can find  $h \in \mathsf{R}_{Z_V}$
with $||h||_{C^k} \leq \eps$.
With this choice of $h$
$$||V - V[h]||_{C^k} \leq \eps$$
and
$X_{V[h]}$ has nondegenerate equilibria. This concludes the proof of the density.

To see that $DZ_V(f)$ is Fredholm,
write $DZ_V(f) = A \circ B \circ C$
where $C : {\cal B}_0^k \to C^k(M), B : C^k(M) \to C^k(M)$ and $A :C^k(M) \to {\cal B}_0^k$ are respectively defined by
$Ch = f.(Vh) + h, \; Bh = \frac{1}{f} h$ and  $A h = h - <h, \one>.$

The operator $C$ is the sum of a compact operator  and identity. Hence, by a classical result, (see e.g Lang, 1993, Theorem 2.1, Chapter XVII) it is Fredholm. Operators  $B$  and $A$ are clearly Fredholm since $Ker(B) = \{0\}, Im(B) = C^k(M), Ker(A) = \RR \one$
and $Im(A) = {\cal B}_0^k.$
Since, the composition of Fredholm operators is Fredholm (Lang, 1993, Corollary 2.6  Chapter XVII), $DZ_V(f)$ is Fredholm. 
\qed
\providecommand{\bysame}{\leavevmode\hbox to3em{\hrulefill}\thinspace}

\end{document}